\documentclass[12pt]{article}
\SetMathAlphabet\mathcal{normal}{U}{rsfs}{m}{n}
\usepackage{a4wide}
\usepackage{amsmath}
\usepackage{theorem}
\usepackage{amsfonts}
\usepackage[english,french]{babel}
\newlength{\secskip}
\newlength{\thsecskip}
\newlength{\ssecskip}
\newlength{\thssecskip}
\newlength{\sssecskip}
\newlength{\thsssecskip}
\makeatletter
\renewcommand\section{\@startsection{section}{1}{\z@}%
  {\secskip}
{2.3ex \@plus.2ex}
{\normalfont\Large\bfseries}}
\renewcommand\subsection{\setcounter{equation}{0}%
\renewcommand{\theequation}{\thesubsection.\arabic{equation}}%
  \@startsection{subsection}{2}{\z@}%
{\ssecskip}
  {1.5ex \@plus .2ex}
{\normalfont\normalsize\bfseries}}
\renewcommand\subsubsection{\setcounter{equation}{0}%
\renewcommand{\theequation}{\thesubsubsection.\arabic{equation}}%
\@startsection{subsubsection}{3}{\z@}%
  {\sssecskip}
  {\z@}
{\normalfont\normalsize\bfseries}}
\makeatother
\newcommand\Subsection[1]{\subsection{{\hskip-0.5em}#1}}
\newcommand\Subsubsection[1]{\subsubsection{{\hskip-0.5em}#1}}
\newcommand\rref[1]{{\rm\ref{#1}}}
\theoremstyle{change}
\theorembodyfont{\sl}
\newtheorem{xthm}[subsection]{Th\'{e}or\`{e}me.}
\newtheorem{xprop}[subsection]{Proposition.}
\newtheorem{xlem}[subsection]{Lemme.}
\newtheorem{xcor}[subsection]{Corollaire.}
\newtheorem{xsublem}[subsubsection]{Lemme.}
\newtheorem{xsubthm}[subsubsection]{Th\'{e}or\`{e}me.}
\newtheorem{xsubprop}[subsubsection]{Proposition.}
\newtheorem{xsubcor}[subsubsection]{Corollaire.}
\newtheorem{xdefi}[subsection]{D\'{e}finition.}
\newtheorem{xsubdefi}[subsubsection]{D\'{e}finition.}
\newtheorem{xsubsublem}{Lemme.}[subsubsection]
\newtheorem{xsubsubprop}{Proposition.}[subsubsection]
\newtheorem{xsubsubcor}{Corollaire.}[subsubsection]
\theorembodyfont{\rmfamily}
\newtheorem{xrem}[subsection]{Remarque.}
\newtheorem{xexam}[subsection]{Exemple.}
\newtheorem{xsubrem}[subsubsection]{Remarque.}
\newtheorem{xsubexam}[subsubsection]{Exemple.}
\newenvironment{thm}{\addvspace{\thssecskip}\begin{xthm}}{\end{xthm}}
\newenvironment{prop}{\addvspace{\thssecskip}\begin{xprop}}{\end{xprop}}
\newenvironment{lem}{\addvspace{\thssecskip}\begin{xlem}}{\end{xlem}}
\newenvironment{cor}{\addvspace{\thssecskip}\begin{xcor}}{\end{xcor}}
\newenvironment{defi}{\addvspace{\thssecskip}\begin{xdefi}}{\end{xdefi}}
\newenvironment{rem}{\addvspace{\thssecskip}\begin{xrem}}{\end{xrem}}

\newenvironment{sublem}{\addvspace{\thsssecskip}\begin{xsublem}}{\end{xsublem}}

\newenvironment{subprop}{\addvspace{\thsssecskip}\begin{xsubprop}}{\end{xsubprop}}
\newenvironment{subcor}{\addvspace{\thsssecskip}\begin{xsubcor}}{\end{xsubcor}}
\newenvironment{subdefi}{\addvspace{\thsssecskip}\begin{xsubdefi}}{\end{xsubdefi}}
\newenvironment{subrem}{\addvspace{\thsssecskip}\begin{xsubrem}}{\end{xsubrem}}

{\end{xsubsublem}}
{\end{xsubsubprop}}
{\end{xsubsubcor}}
\newcounter{nc}
\renewcommand{\thenc}{{\rm(\roman{nc})}}
\newenvironment{romlist}%
{\begin{list}{\thenc}{
\usecounter{nc}
\parsep=0pt
\setlength  \labelwidth{\leftmargin}
\addtolength\labelwidth{-\labelsep}
}
}{\end{list}}
\newcounter{nnc}
\renewcommand{\thennc}{{\rm(\alph{nnc})}}
\newenvironment{subromlist}%
{\begin{list}{\thennc}{
\usecounter{nnc}
\parsep=0pt
\setlength  \labelwidth{\leftmargin}
\addtolength\labelwidth{-\labelsep}
}
}{\end{list}}
{\begin{list}{$\bullet$}{\parsep=0pt}}{\end{list}}
\renewcommand{\ker}{{\rm Ker\:}}
\newcommand{\car}{{\rm car\:}}
\newcommand{\Spec}{{\rm Spec\:}}
\newcommand{\NN}{\mathbb{N}}
\newcommand{\ZZ}{\mathbb{Z}}
\newcommand{\QQ}{\mathbb{Q}}
\newcommand{\CC}{\mathbb{C}}
\newcommand{\Aa}{\mathbb{A}}

\newcommand{\PP}{\mathbb{P}}
\newcommand{\GG}{\mathbb{G}}
\newcommand{\soul}[1]{\underline{#1}}
\newcommand{\gal}{{\rm Gal}}
\newcommand{\gals}[1]{{\rm Gal\,}(#1_{{\rm s}}/#1)}

\newcommand{\hypG}[1]{(2)}
\newcommand{\cyc}[1]{\textup{Cyc($#1$)}}
\newcommand{\HH}{\mathrm{H}}
\newcommand{\re}{\mathrm{R}}
\newcommand{\ree}{\mathrm{R_{el}}}
\newcommand{\rre}{\overset{\re}{\sim}}
\newcommand{\rree}{\overset{\ree}{\sim}}
\newcommand{\gre}{\mathfrak{R}}
\newcommand{\gree}{\mathfrak{R}_\mathrm{el}}
\newcommand{\ssoul}[1]%
{\rlap{$\underline{#1}$}%
\underline{\phantom{\lower.25ex\hbox{$#1$}}}}

\newcommand{\carrenoir}{\vrule height2.5mm width2mm depth 0mm}
\def\qed{ \hbox to5mm{}\hfill\llap{\carrenoir}}

\newcommand{\dem}{{\sl \noindent D\'{e}monstration. }}
\newcommand{\fl}{\rightarrow}
\newcommand{\ffl}{\longrightarrow}
\newcommand{\inj}{\hookrightarrow}

\newcommand{\sfl}[1]{\mathop{\fl}\limits^{#1}}
\newcommand{\sffl}[1]{\mathop{\ffl}\limits^{#1}}
\newcommand{\sinj}[1]{\mathop{\inj}\limits^{#1}}
\newcommand{\flis}{\sfl{\sim}}

\newcommand{\surj}{\rlap{$\rightarrow$}\kern.1pt\rightarrow}
\newcommand{\ssurj}{\rlap{$\longrightarrow$}\kern.1pt\longrightarrow}
\newcommand{\surjgauche}{\rlap{$\leftarrow$}\kern.1pt\leftarrow}
\newcommand{\ssurjgauche}{\rlap{$\longleftarrow$}\kern.1pt\longleftarrow}
\newcommand{\mapdown}[1]%
{\big\downarrow\rlap{$\vcenter{\hbox{$\scriptstyle#1$}}$}}
\newcommand{\longmapstodown}%
{\vbox{\hbox to 3mm{\hss$\scriptscriptstyle\top$\hss}\hbox to
3mm{\hss$\big\downarrow$\hss}}}
\newcommand{\varfl}[1]{%
\setbox0=\hbox{$\;\;{\scriptstyle#1}\;\;\;$}%
\setbox1=\hbox to\wd0{$\;$\rightarrowfill$\;$}%
{\mathop{\box1}\limits^{\box0}}%
}
\newcommand{\varflspecial}[2]{%
\setbox0=\hbox{$\;\;{\scriptstyle#1}\;\;\;$}
\setbox1=\hbox to\wd0{$\;$\rightarrowfill$\;$}
{\mathop{\box1}\limits^{#2}}
}
\newcommand{\varflman}[2]{%
\setbox0= \hbox to #1{$\;$\rightarrowfill$\;$}
{\mathop{\box0}\limits^{#2}}
}
\newcommand{\varmapstoman}[2]{%
\setbox0= \hbox to #1{$\;\mapstochar$\rightarrowfill$\;$}
{\mathop{\box0}\limits^{#2}}
}

\date{}
\title{Sur la R-\'{e}quivalence de torseurs sous un groupe fini
\thanks{Accept\'{e} pour publication au
\textsl{Journal of Number Theory\/} (20 septembre 2002).}}
\author{Laurent Moret-Bailly
\thanks{L'auteur est membre du r\'{e}seau europ\'{e}en \og Arithmetic
Algebraic Geometry\fg\ (contrat HPRN-CT-2000-00120).
}\medskip
\\
{\small IRMAR (Institut de Recherche Math\'{e}matique de Rennes, UMR
6625 du CNRS)}\\
{\small Universit\'{e} de Rennes 1, Campus de Beaulieu, F-35042 Rennes Cedex}\\
{\small Laurent.Moret-Bailly@univ-rennes1.fr}\\
{\small http://www.maths.univ-rennes1.fr/\lower1ex\hbox{\~{}}moret/}
}
\begin{document}
\maketitle
\selectlanguage{french}
%
\begin{abstract}\medskip

On donne des crit\`{e}res de R-\'{e}quivalence pour les torseurs sous un
sch\'{e}ma en groupes fini constant sur un corps. En particulier, un
\'{e}nonc\'{e} de d\'{e}vissage galoisien permet, en utilisant la notion
de bitorseur, de formaliser et de g\'{e}n\'{e}raliser un th\'{e}or\`{e}me
de Philippe Gille dans le cas des corps locaux: ce dernier est notamment
\'{e}tendu aux corps locaux sup\'{e}rieurs.
\end{abstract}
\selectlanguage{english}
\begin{abstract}
We give criteria for R-equivalence of torsors under finite constant
group schemes over a field. In particular, using bitorsors, we obtain a Galois
d\'{e}vissage result which formalises and generalises a theorem of
Philippe Gille in the case of local fields; for instance, Gille's theorem
is shown to extend to higher local fields.
\end{abstract}
\selectlanguage{french}

\noindent\textsl{Classification AMS 2000:}
12G05, 
14G20, 
18G50. 

\section{Introduction}\label{Introd}
\Subsection{R-\'{e}quivalence.}\label{RapREquiv} Soient $K$ un corps
et $A$ l'anneau semi-local de la droite affine $\Aa^1_{K}$ en
$\{0,1\}$. On rappelle
que si $H$ est un $K$-sch\'{e}ma en groupes, deux
$H$-torseurs $X$ et $Y$ sur $K$ ({\it i.e.\/} sur $\Spec K$, et pour
la topologie fppf) sont dits \emph{\'{e}l\'{e}mentairement
R-\'{e}quivalents\/} s'il existe
un $H$-torseur sur $\Spec A$ induisant (\`{a} isomorphisme pr\`{e}s)
$X$ en $0$ et $Y$ en $1$. La R-\'{e}quivalence est par d\'{e}finition
la relation
d'\'{e}quivalence engendr\'{e}e par la R-\'{e}quivalence
\'{e}l\'{e}mentaire, sur la classe des
$H$-torseurs sur $K$. Elle induit naturellement une relation
d'\'{e}quivalence, encore appel\'{e}e R-\'{e}quivalence, sur l'ensemble
$\HH^1(K,H)$ des \emph{classes d'isomorphie} de $H$-torseurs sur $K$
(le lien avec la d\'{e}finition de $\HH^1(K,H)$ par cocycles est \'{e}tabli,
avec toute la g\'{e}n\'{e}ralit\'{e} voulue, dans \cite{giraud}, III,
3.6.4 et 3.6.5,
ainsi, pour la cohomologie galoisienne, que dans \cite{D-G}, III,
\S\,5, n$^{\rm os}$~3 et~4).

On a en particulier des sous-ensembles
$$
\ree(K,H)\subset\re(K,H)\subset\HH^1(K,H)
$$
o\`{u} $\re(K,H)$ (resp. $\ree(K,H)$) est l'ensemble des classes de
$H$-torseurs
R-\'{e}quivalents (resp. \'{e}l\'{e}mentairement R-\'{e}quivalents)
au torseur trivial.

D'autre part, si $L$ est une extension de $K$, on note comme
d'habitude
$\mathrm{H}^1(L/K,H)\subset\mathrm{H}^1(K,H)$  l'ensemble des classes
de
$H$-torseurs trivialis\'{e}s par $L$.

Dans ce qui suit, on consid\'{e}rera surtout des sch\'{e}mas en
groupes
finis \'{e}tales sur $K$, et des torseurs sous iceux; si l'on fixe
une
cl\^{o}ture s\'{e}parable $K_{\rm s}$ de $K$ et que l'on pose
$\gal_{K}=\gals{K}$, on peut identifier un $K$-sch\'{e}ma fini
\'{e}tale $X$
au $\gal_{K}$-ensemble fini $X(K_{\rm s})$, ce que nous ferons
syst\'{e}matiquement. Noter que $\gal_{K}$ lui-m\^{e}me, ainsi
que ses sous-groupes distingu\'{e}s et ses quotients, peut \^{e}tre
consid\'{e}r\'{e}
comme un $K$-sch\'{e}ma en groupes profini, lorsqu'on le munit de son
action sur lui-m\^{e}me par conjugaison.

\Subsection{Notations.}\label{NotThPpal} Soient $K$, $K_{\rm s}$ et
$\gal_{K}$ comme ci-dessus, et soit $M\subset K_{\rm s}$ une
extension galoisienne de $K$, non n\'{e}cessairement finie, de groupe
$\Pi$ (qui est donc un groupe profini, quotient de $\gal_{K}$ par
un
sous-groupe ferm\'{e} distingu\'{e}). On suppose donn\'{e}e une suite exacte
\begin{equation}\label{DevPi}
1\ffl \Gamma \ffl \Pi \ffl \pi \ffl 1
\end{equation}
de groupes profinis (ici et dans la suite, tous les morphismes
de groupes profinis sont suppos\'{e}s continus); on note $M_0\subset M$
le corps des invariants
de $\Gamma$, de sorte que $\pi=\mathrm{Gal}\,(M_0/K)$.

Enfin on se donne un groupe fini $G$, que l'on voit comme
$K$-sch\'{e}ma en groupes constant.

\begin{thm} \label{ThPpal}
Avec les notations de {\rm\ref{NotThPpal}}, on suppose que:
\begin{romlist}
\item\label{ThPpal0} la suite exacte \textup{(\ref{DevPi})} de
groupes profinis est \emph{scind\'{e}e};
\item\label{ThPpal1} on a $\HH^1(M_0/K,G)\subset\re(K,G)$;
\item\label{ThPpal2} pour tout $K$-sch\'{e}ma en groupes $H$ quotient
de $\Gamma$, dont le groupe sous-jacent est isomorphe \`{a} un
sous-groupe de $G$, on a
$\HH^1(M/K,H)\subset\re(K,H)$.
\end{romlist}
Alors on a $\HH^1(M/K,G)\subset\re(K,G)$.
\end{thm}

\Subsection{Remarques sur les hypoth\`{e}ses. }\label{RemThPpal}

\Subsubsection{}\label{RemThPpal1} Comme le lecteur pourra le constater
(voir la preuve du th\'{e}or\`{e}me \ref{ThDeviss} et plus
pr\'{e}cis\'{e}ment le diagramme (\ref{DevG'})) on peut remplacer
l'hypoth\`{e}se \ref{ThPpal0} par la suivante, plus faible: pour tout
morphisme continu $\theta:\Pi\to G$, le morphisme naturel
$\overline{\theta}:\pi\to \theta(\Pi)/\theta(\Gamma)$ (d\'{e}duit de
$\theta$ par passage au quotient) se rel\`{e}ve en un morphisme
continu $\pi\to\theta(\Pi)$.

\Subsubsection{}\label{RemThPpal2}
La preuve utilise de fa\c{c}on cruciale l'identification
de $\HH^1(M/K,G)$ (pour $G$ constant) avec le quotient de
$\mathrm{Hom}\,(\Pi,G)$ par la conjugaison. De ce fait, elle ne
s'\'{e}tend pas aux sch\'{e}mas en groupes
\'{e}tales g\'{e}n\'{e}raux, \`{a} l'exception, de mani\`{e}re assez
formelle, de certaines \emph{formes int\'{e}rieures} de groupes
constants. Plus pr\'{e}cis\'{e}ment,
avec les hypoth\`{e}ses de \ref{ThPpal}, soit de plus $X$ un $G$-torseur
(\`{a} droite, disons) \emph{trivialis\'{e} par $M$} et soit $G'$ le
$K$-sch\'{e}ma en groupes $\soul{\mathrm{Aut}}_{G}(X)$:
alors on a  une bijection de $\HH^1(K,G)$ avec $\HH^1(K,G')$,
associant au $G$-torseur \`{a} droite $Y$ le $K$-sch\'{e}ma fini
\'{e}tale $\soul{\mathrm{Isom}}_{G}(X,Y)$, muni de l'action \`{a} droite
\'{e}vidente de $G'$. On voit facilement que cette bijection respecte
la R-\'{e}quivalence et envoie $\HH^1(M/K,G)$ sur $\HH^1(M/K,G')$, de sorte
que l'on d\'{e}duit encore de \ref{ThPpal} que 
$\HH^1(M/K,G')\subset\re(K,G')$.

\Subsubsection{}\label{RemThPpal3}
On peut pr\'{e}ciser l'\'{e}nonc\'{e} en tenant compte
de la filtration naturelle sur l'ensemble point\'{e} $\HH^1(K,H)$.
Pour $n\in\NN$,
d\'{e}finissons $\re_{n}(K,H)\subset\HH^1(K,H)$ comme suit:
$\re_{0}(K,H)$ est r\'{e}duit \`{a} la classe triviale, et
$\re_{n+1}(K,H)$ est l'ensemble des classes de $H$-torseurs
\'{e}l\'{e}mentairement R-\'{e}quivalents \`{a} un torseur de
$\re_{n}(K,H)$. Si l'on remplace l'inclusion de
\ref{ThPpal}\,\ref{ThPpal1} par $\HH^1(M_0/K,G)\subset\re_{m}(K,G)$
et celle de \ref{ThPpal}\,\ref{ThPpal2} par
$\HH^1(M/K,H)\subset\re_{n}(K,H)$, alors on peut conclure que
$\HH^1(M/K,G)\subset\re_{m+n}(K,G)$. Ceci r\'{e}sulte facilement du
th\'{e}o\-r\`{e}me de d\'{e}vissage \ref{ThDeviss} et de la remarque
\ref{RemOrdreREquiv}.

\Subsubsection{}\label{RemThPpal4}
Dans la situation de \ref{ThPpal}, on pourrait envisager
une relation  plus fine que la  R-\'{e}qui\-va\-lence, \`{a} savoir la
\og R$_{M/K}$-\'{e}quivalence\fg, d\'{e}finie comme la relation
d'\'{e}quivalence \emph{dans $\HH^1(M/K,G)$} engendr\'{e}e par la
R-\'{e}quivalence \'{e}l\'{e}mentaire. Si l'on note $\re(M/K,G)$
l'ensemble des $G$-torseurs R$_{M/K}$-\'{e}quivalents au torseur
trivial, le lecteur pourra constater que
si l'on remplace dans \ref{ThPpal} l'hypoth\`{e}se \ref{ThPpal1} par
$\HH^1(M_0/K,G)\subset\re(M/K,G)$, et \ref{ThPpal2} par
$\HH^1(M/K,H)\subset\re(M/K,H)$, alors on conclut que
$\HH^1(M/K,G)\subset\re(M/K,G)$.  On peut de plus combiner
cette variante de \ref{ThPpal} avec celle de \ref{RemThPpal3} ci-dessus.

\Subsubsection{}\label{RemThPpal5}
On n'a pas inclus dans l'\'{e}nonc\'{e} les remarques
\ref{RemThPpal3} et \ref{RemThPpal4} qui pr\'{e}c\`{e}dent, pour
la raison suivante: dans les applications les plus importantes,
$K$ est un corps valu\'{e} hens\'{e}lien, donc est fertile (pour 
toute $K$-vari\'{e}t\'{e} $V$ lisse connexe, $V(K)$ est vide ou dense 
dans $V$), et dans ce cas
le th\'{e}or\`{e}me 2 de \cite{ReqSim} montre que la R-\'{e}quivalence
dans $\HH^1(K,G)$ co\"{\i}ncide avec la R-\'{e}quivalence
\'{e}l\'{e}mentaire, ce qui  entra\^{\i}ne en outre que sa restriction
\`{a} $\HH^1(M/K,G)$ co\"{\i}ncide avec la R$_{M/K}$-\'{e}quivalence. 
On voit donc
que \ref{RemThPpal3} et \ref{RemThPpal4} sont sans objet
dans ce cas.
\medskip

Le th\'{e}or\`{e}me \ref{ThPpal} sera \'{e}tabli au \S\,\ref{SecDem};
auparavant, au \S\,\ref{SecMod}, nous en d\'{e}duirons le suivant:

\begin{thm}\label{ThRec} Soit $K$ le corps des fractions d'un anneau
	de valuation discr\`{e}te hens\'{e}lien, \`{a} corps r\'{e}siduel
	$k$ de caract\'{e}ristique $p\geq0$, et soit $G$ un groupe fini
	\emph{d'ordre non divisible par $p$}.

	On suppose que $\re(k,G)=\HH^1(k,G)$. Alors $\re(K,G)=\HH^1(K,G)$
	\emph{(et donc $\ree(K,G)=\HH^1(K,G)$ d'apr\`{e}s
	\cite{ReqSim}, cf. remarque \ref{RemThPpal5})}.
\end{thm}

\begin{subrem}\label{RemCorpsResFini}  Bien entendu, $G$ est
consid\'{e}r\'{e}
	ici comme groupe constant sur $k$ et sur $K$ respectivement.

	L'hypoth\`{e}se de l'\'{e}nonc\'{e} est notamment v\'{e}rifi\'{e}e
	lorsque $k$ est \emph{fini}, comme
	il r\'{e}sulte par exemple
	de \ref{CorGalProCyc} plus bas. On retrouve ainsi le th\'{e}or\`{e}me 1
	de \cite{gille}, qui a servi de point de d\'{e}part au pr\'{e}sent
	travail.

	Mais l'hypoth\`{e}se est aussi v\'{e}rifi\'{e}e, 
trivialement, lorsque $k$ est
	\emph{s\'{e}parablement clos}, et plus g\'{e}n\'{e}ralement 
(\ref{CorGalProCyc})
	lorsque $\gal_{k}$ est pro-cyclique.

\end{subrem}

\begin{subrem}\label{ContrExZ/8} L'\'{e}nonc\'{e} serait faux sans
	l'hypoth\`{e}se
	sur l'ordre de $G$, comme le montre l'exemple suivant de
	Colliot-Th\'{e}l\`{e}ne (\cite{jlct1}, Appendix): soient $K=\QQ_{2}$,
	$G=\ZZ/8\ZZ$, et $L$ l'unique extension non ramifi\'{e}e de degr\'{e}
	$8$ de $K$. Alors
	$\Spec L$ est un $G$-torseur sur $K$, qui n'est pas R-\'{e}quivalent
	au torseur trivial, alors que tout $G$-torseur sur le corps 
r\'{e}siduel
	de $K$ l'est (cf. \ref{CritTVR}\,\ref{CritTVR1} plus loin).
\end{subrem}

Par r\'{e}currence, on en d\'{e}duit  un r\'{e}sultat analogue pour
les \og corps locaux
sup\'{e}rieurs\fg; plus pr\'{e}cis\'{e}ment:

\begin{subcor}\label{CorpsLocSup} Consid\'{e}rons une suite de corps
	$$K=K_{0},\;K_{1},\ldots,\;K_{n}$$
	o\`{u}, pour chaque $i\in\{0,\ldots,n-1\}$, $K_{i}$ est muni d'une
	valuation discr\`{e}te hens\'{e}lienne de corps r\'{e}siduel
	$K_{i+1}$,
	et o\`{u} le groupe de Galois absolu de $K_{n}$ est pro-cyclique
	\emph{(condition v\'{e}rifi\'{e}e notamment si $K_{n}$ est fini ou
	s\'{e}parablement clos)}.

	Alors, si $G$ est un groupe fini d'ordre premier \`{a} la
	caract\'{e}ristique de $K_{n}$, on a $\re(K,G)=\HH^1(K,G)$
	\emph{(et donc encore $\ree(K,G)=\HH^1(K,G)$, d'apr\`{e}s
	\cite{ReqSim} si $n>0$ et d'apr\`{e}s \ref{CorGalProCyc} si
	$n=0$)}.
\qed
\end{subcor}

\begin{subrem}\label{RemSerFormelles} Par exemple, on a
$\re(K,G)=\HH^1(K,G)$ pour tout groupe fini $G$, lorsque $K$ est un
corps de s\'{e}ries formelles de la forme 
$\CC(\!(X_{1})\!)\cdots(\!(X_{n})\!)$;
m\^{e}me ce cas particulier semble nouveau.
\end{subrem}

\Subsection{Plan de l'article. }\label{plan}
Au \S\,\ref{SecREquiv}, nous donnons des crit\`{e}res simples (et
sans doute plus ou moins bien connus) de R-\'{e}quivalence.

Le \S\,\ref{SecMod} est consacr\'{e} \`{a} la preuve du
th\'{e}or\`{e}me \ref{ThRec} \`{a} partir de \ref{ThPpal}.

Au
\S\,\ref{SecResAb} on donne d'autres applications de \ref{ThPpal},
lorsque $K$ est valu\'{e} hen\-s\'{e}\-lien \`{a} corps
r\'{e}\-si\-duel $k$ ab\'{e}\-lien. Ainsi, si $k$ est fini de 
caract\'{e}ristique
dif\-f\'{e}\-rente de $2$, et si $G$ est un groupe fini quelconque, on montre
que, avec les notations tra\-di\-tion\-nelles,  $\HH^1(K_{\rm mod}K_{\rm
ab}/K,G)\subset\ree(K,G)$ (voir \ref{ThResAb} pour un \'{e}nonc\'{e} 
plus pr\'{e}cis).

La preuve de \ref{ThPpal} occupe le \S\,\ref{SecDem}; elle
s'inspire de \cite{gille} et repose sur un d\'{e}vissage qui est
d\'{e}crit
dans loc. cit. en termes de cocycles, et que l'on traduit ici, de
mani\`{e}re plus fonctorielle, dans le langage des bitorseurs. Comme
ceux-ci n'ont pas encore la m\^{e}me popularit\'{e} que les torseurs,
on pr\'esente divers sorites sur cette notion
(\S\,\ref{SecBitors}), son comportement vis-\`{a}-vis de la
R-\'{e}quivalence (\S\,\ref{SecBitorsREquiv}) et le cas o\`{u} l'un
des groupes est constant (\S\,\ref{SecBitorsGal}).

\section{R-\'{e}quivalence: quelques crit\`{e}res
simples}\label{SecREquiv}
On donne ici quelques conditions suffisantes simples de
R-\'{e}quivalence \'{e}l\'{e}mentaire pour les torseurs sous un
sch\'{e}ma en groupes, le plus souvent fini, sur un corps.

Dans tout ce paragraphe, $K$ d\'{e}signe un corps, $K_{{\rm s}}$ une
cl\^{o}ture s\'{e}parable de $K$, et $\gal_{K}=\gals{K}$. Tous les
sch\'{e}mas en groupes consid\'{e}r\'{e}s seront affines de type
fini, et
\og $G$-torseur\fg\ signifie \og $G$-torseur pour la topologie
fppf\fg.
On note ${\rm cd}(K)$ la dimension cohomologique de $K$.
\smallskip

Plut\^{o}t que la R-\'{e}quivalence, nous consid\'{e}rerons ici la
propri\'{e}t\'{e}
suivante:

\begin{defi}\label{DefTVR} Soient $k$ un corps et $G$ un $k$-sch\'{e}ma
en groupes.
	On dit
	que $G$ v\'{e}rifie la propri\'{e}t\'{e} {\rm(TVR)} (\og torseur versel
	rationnel\fg) s'il existe un ouvert $U$ d'un espace affine sur $k$ et
	un $G\times_{k}U$-torseur $\pi:V\to U$ tels que l'application 
naturelle de
	$U(k)$ dans ${\rm H}^1(k,G)$ d\'{e}duite de $\pi$ soit
\emph{surjective}.
\end{defi}

\begin{subrem} La  propri\'{e}t\'{e} {\rm(TVR)} entra\^{\i}ne
	imm\'{e}diatement que $\ree(k,G)=\HH^1(k,G)$. Lorsque $k$ est
	\emph{infini}, elle donne un peu
	mieux: \'{e}tant donn\'{e}e
	une famille finie $(X_{1},\ldots,X_{r})$ de $G$-torseurs, il existe
	un ouvert $U$ de $\Aa^1_{k}$, un $G$-torseur $f:V\to U$ et des
	points rationnels $u_{1},\ldots,u_{n}\in U(k)$ tels que
$f^{-1}(u_{i})\cong
	X_{i}$ pour tout $i$.
\end{subrem}

\begin{subrem} Il est naturel de consid\'{e}rer des variantes de
{\rm(TVR)}, notamment celle obtenue en rempla\c{c}ant \og ouvert
d'espace affine\fg\ par \og vari\'{e}t\'{e} rationnelle\fg. Les
con\-si\-d\'{e}\-ra\-tions ci-dessous resteraient valables, mais
l'\'{e}nonc\'{e} \ref{CritTVR}
serait plus faible.
\end{subrem}

\begin{prop}\label{SorTVR} Soient $G_{i}$ ($1\leq i\leq3$) des
	$K$-sch\'{e}mas en groupes.
	\begin{romlist}
		\item\label{SorTVR1}  Si $G_{1}$ et $G_{2}$ v\'{e}rifient
{\rm(TVR)}, il en est de m\^{e}me de $G_{1}\times G_{2}$.
		\item\label{SorTVR2}  Soit $\varphi:G_{1}\to G_{2}$ 
un morphisme de
		$K$-sch\'{e}mas en groupes. On suppose que l'application
		${\rm H}^1(K,\varphi):
		{\rm H}^1(K,G_{1})\to {\rm H}^1(K,G_{2})$ est 
surjective, et que
$G_{1}$
		v\'{e}rifie {\rm(TVR)}. Alors $G_{2}$ v\'{e}rifie {\rm(TVR)}.
		\item\label{SorTVR3} On suppose que $G_{1}$ est un sous-groupe
		ferm\'{e} de $G_{2}$, et que:
		\begin{subromlist}
			\item\label{SorTVR31} l'inclusion $G_{1}\inj 
G_{2}$ induit
			l'application \emph{triviale}
			${\rm H}^1(K,G_{1})\to{\rm H}^1(K,G_{2})$;
			\item\label{SorTVR32} le $K$-sch\'{e}ma 
$G_{2}/G_{1}$ est un
ouvert
			d'espace affine.
		\end{subromlist}
		Alors  $G_{1}$ v\'{e}rifie {\rm(TVR)}.
		\item\label{SorTVR4} Soit $L$ une extension finie de 
$K$, et soit
		$H$ un $L$-sch\'{e}ma en groupes v\'{e}rifiant 
{\rm(TVR)}. Alors la
		restriction de Weil ${\rm R}_{L/K}(H)$ v\'{e}rifie {\rm(TVR)}.
	\end{romlist}
\end{prop}
\dem Les assertions \ref{SorTVR1} et \ref{SorTVR2}
sont faciles (et la r\'{e}ciproque de \ref{SorTVR1} est un cas particulier
de \ref{SorTVR2}).

Montrons \ref{SorTVR3}. Posons $V=G_{2}$, $U=G_{2}/G_{1}$ (qui est
un ouvert d'espace affine par l'hypoth\`{e}se \ref{SorTVR32}),
et soit $\pi:V\to
U$ le morphisme canonique, qui fait de $V$ un
$G_{1}\times_{k}U$-torseur \`{a} droite pour l'action naturelle de
$G_{1}$ sur $G_{2}$. Soit $X$ un $G_{1}$-torseur (\`{a} droite) sur $K$:
l'hypoth\`{e}se \ref{SorTVR31} montre que $X$ se plonge dans le
$G_{2}$-torseur trivial, de mani\`{e}re compatible aux actions de $G_{1}$;
ceci \'{e}quivaut \`{a} dire que $X$ (ainsi plong\'{e}) est une fibre 
de $\pi$ en
un point de $U(K)$.

Pour \ref{SorTVR4}, soit $V\to U$ un $H$-torseur
comme dans la d\'{e}finition \ref{DefTVR}, o\`{u} $U$ est un ouvert
d'espace affine sur $L$: alors on en d\'{e}duit par restriction de Weil un
${\rm R}_{L/K}(H)$-torseur $V_{1}\to U_{1}$, o\`{u} $U_{1}$ est un ouvert
d'espace affine sur $K$, dont on v\'{e}rifie qu'il a la propri\'{e}t\'{e}
voulue.\qed

\Subsection{La condition {\mdseries\cyc{K,G}}.}\label{DefCondCyc}
\begin{subdefi}\label{Def2Exp} Soit $G$ un groupe fini. Le
\emph{$2$-exposant} de $G$ est l'ordre maximum d'un \'{e}l\'{e}ment
de $2$-torsion de $G$.
\end{subdefi}

Si $G$ est un groupe fini et $K$ un corps, nous aurons \`{a}
consid\'{e}rer la condition suivante (o\`{u} l'on convient que
$K(\mu_{2^e})=K$ si $\car K=2$):
\bigskip

\noindent\cyc{K,G}:\hskip2em si $2^e$ d\'{e}signe le $2$-exposant de
$G$, l'extension $K(\mu_{2^e})/K$ est cyclique.

\begin{subrem}\label{RemCPCyc} La condition \cyc{K,G} est
v\'{e}rifi\'{e}e
dans chacun des cas suivants:
\begin{romlist}
\item\label{RemCPCyc1} $\car K>0$;
\item\label{RemCPCyc2} $-1$ ou $-2$ est un carr\'{e} dans $K$;
\item\label{RemCPCyc3} $K$ est le corps des fractions d'un anneau
local int\`{e}gre hens\'{e}lien \`{a} corps r\'{e}siduel $k$ de
caract\'{e}ristique
diff\'{e}rente de $2$, tel que \cyc{k,G} soit v\'{e}rifi\'{e}e;
\item\label{RemCPCyc4} $G$ n'a pas d'\'{e}l\'{e}ment d'ordre $8$.
\end{romlist}
\smallskip

Remarquer d'autre part que si $l$ est un nombre premier \emph{impair,
toute} extension de la forme $K(\mu_{l^m})/K$ est cyclique; nous
utiliserons sans commentaire cette propri\'{e}t\'{e}.
\end{subrem}

\begin{thm}\label{CritTVR} Soit $G$ un $K$-sch\'{e}ma en groupes
	\emph{commutatif}. Dans
	chacun des cas sui\-vants, $G$ v\'{e}rifie {\rm(TVR)}:
	\begin{romlist}
		\item\label{CritTVR1} $K$ est de caract\'{e}ristique 
$p>0$, et $G$
		est un $p$-groupe fini constant;
		\item\label{CritTVR2} $G$ est de type multiplicatif 
d\'{e}ploy\'{e}
		par une extension m\'{e}tacyclique de $K$ \emph{(on
		rappelle qu'un groupe fini est dit m\'{e}tacyclique si ses
		sous-groupes de Sylow sont cycliques)};
		\item\label{CritTVR3} $G$ est de type multiplicatif et
		${\rm cd}(K)\leq1$;
		\item\label{CritTVR4} $G$ est fini \'{e}tale  et ${\rm cd}(K)\leq1$;
		\item\label{CritTVR5} $G$ est fini constant, et la condition
		\cyc{K,G} de \rref{DefCondCyc} est satisfaite.
	\end{romlist}
\end{thm}
\dem Rappelons qu'un $K$-sch\'{e}ma en groupes $G$ est \emph{de type
multiplicatif d\'{e}ploy\'{e}} s'il est isomorphe \`{a} un sous-sch\'{e}ma
en groupes
ferm\'{e} de $\GG_{{\rm m},K}^n$, pour $n$ convenable; il est \emph{de type
multiplicatif} si $G_{K_{{\rm s}}}$ est de type
multiplicatif d\'{e}ploy\'{e} comme $K_{{\rm s}}$-sch\'{e}ma en groupes.
\smallskip

\ref{CritTVR1} D'apr\`{e}s \ref{SorTVR}\,\ref{SorTVR1}, on peut
supposer que $G=\ZZ/p^n\ZZ$. Soit alors $W$ le $K$-sch\'{e}ma en groupes
des vecteurs de Witt tronqu\'{e}s de longueur $n$. On a une suite exacte
\og d'Artin-Schreier-Witt\fg
$$0 \ffl G \varfl{\phantom{\Phi-{\rm Id}}} W \varfl{\Phi-{\rm Id}} W
\ffl 0$$
o\`{u} $\Phi$ est l'endomorphisme de Frobenius. Comme $W$ est
extension
successive de groupes additifs $\GG_{{\rm a},K}$, on a ${\rm
H}^{1}(K,W)=0$; comme $W$ est isomorphe \`{a}
$\Aa^n_{K}$ comme $K$-sch\'{e}ma, on conclut par
\ref{SorTVR}\,\ref{SorTVR3}.
\smallskip

\noindent\ref{CritTVR2} est la proposition 1 de \cite{gille}.
Rappelons
l'argument. Soit $M$ une extension finie de $K$ d\'{e}ployant $G$.
D'apr\`{e}s (\cite{CT-S2}, proposition 1.3), il existe une suite
exacte
$$ 1 \ffl G \ffl S \sffl{\pi} E \ffl 1 $$
dans laquelle:
\begin{itemize}
	\item $E$ est un $K$-tore quasi-trivial, {\it i.e.\/} produit 
de 
	restrictions de
	Weil ${\rm R}_{L_{i}/K}\GG_{{\rm m},L_{i}}$ o\`{u} les $L_{i}$ sont
des
	extensions finies s\'{e}parables de $K$; en particulier $E$ est un
ouvert
	d'espace affine sur $K$;
	\item le tore $S$ (ainsi d'ailleurs que $E$) est d\'{e}ploy\'{e} par
	$M$, et $S$ est \og flasque\fg; si $M/K$ est 
m\'{e}tacyclique, cela implique
	que ${\rm H}^1(K,S)$ est trivial (\cite{CT-S1}, cor. 3).
\end{itemize}
On conclut donc \`{a} nouveau par \ref{SorTVR}\,\ref{SorTVR3}.
\smallskip

\noindent\ref{CritTVR3} L'argument est le m\^{e}me que pour
\ref{CritTVR2}, sans
l'extension m\'{e}tacyclique $M$; ici c'est l'hypoth\`{e}se ${\rm
cd}(K)\leq1$ qui
assure que  ${\rm H}^1(K, S)=1$ (un tore est un groupe divisible).
\smallskip

\noindent\ref{CritTVR4} (signal\'{e} par Philippe Gille):
On \'{e}crit $G$ comme quotient d'un
\og$\gal_{K}$-module de per\-mu\-ta\-tion\fg, c'est-\`{a}-dire d'un produit
$P$
de restrictions de Weil ${\rm R}_{L_{i}/K}(C_{i})$ o\`{u} les $L_{i}$
sont des
extensions finies s\'{e}parables de $K$ et o\`{u}
$C_{i}=(\ZZ/n_{i}\ZZ)_{L_{i}}$.  L'hypoth\`{e}se sur $K$ implique
alors que
l'application naturelle ${\rm H}^1(K,P)\to{\rm H}^1(K,G)$ est
surjective; appliquant successivement les assertions \ref{SorTVR2},
\ref{SorTVR1} et \ref{SorTVR4} de \ref{SorTVR}, on est ramen\'{e} au
cas o\`{u}
$G=\ZZ/l^m\ZZ$ (avec $l$ premier). On conclut par \ref{CritTVR1} si
$l=\car K$, et par \ref{CritTVR3} sinon.
\smallskip

\noindent\ref{CritTVR5} Par d\'{e}composition en produit, on se
ram\`{e}ne encore au cas o\`{u} $G=\ZZ/l^m\ZZ$ ($l$ premier). On
conclut par
\ref{CritTVR1} si $l=\car K$, et par \ref{CritTVR2} sinon ($G$ est
d\'{e}ploy\'{e} par
$K(\mu_{l^m})$).\qed
\medskip

Voici enfin une application aux groupes non n\'{e}cessairement
commutatifs, lorsque $K$ est par exemple un corps fini:

\begin{cor}\label{CorGalProCyc} Soit $G$ un groupe fini, vu comme
$K$-sch\'{e}ma en groupes constant. On suppose que $\gal_K$ est
commutatif, et que la condition \cyc{K,G} de \rref{DefCondCyc} est
satisfaite. Alors $\ree(K,G)=\HH^1(K,G)$.
\end{cor}
\dem Les $G$-torseurs sont classifi\'{e}s par les homomorphismes
$\gal_K\to G$, \`{a} conjugaison pr\`{e}s. Tout $G$-torseur est donc
induit
par un torseur sous un sous-groupe commutatif $H$ de $G$. Comme
\cyc{K,H} est clairement satisfaite, il suffit d'appliquer
\ref{CritTVR}\,\ref{CritTVR5}.\qed

\section{Le cas \og premier \`{a} $p$\fg: preuve du th\'{e}or\`{e}me
\ref{ThRec}}\label{SecMod}
\Subsection{Notations.}\label{NotVal} Soient $K$, $K_{\rm s}$ et
$\gal_{K}$ comme dans \ref{RapREquiv}; on suppose en outre que $K$
est le corps des fractions d'un anneau de valuation discr\`{e}te
hens\'{e}lien $\Lambda$, de corps r\'{e}siduel $k$; on note $p$
l'exposant
caract\'{e}ristique de $k$.

La fermeture int\'{e}grale de $\Lambda$ dans $K_{\rm s}$ est un
anneau de
valuation, dont le corps r\'{e}siduel est une cl\^{o}ture
alg\'{e}brique
$\overline{k}$ de $k$. On note $k_{\rm s}$ la cl\^{o}ture
s\'{e}parable
de $k$ dans $\overline{k}$, et l'on pose $\gal_{k}=\gals{k}$. On a
une suite
exacte canonique

\begin{equation}\label{DevGal}
1\ffl I \ffl \gal_{K} \ffl \gal_{k}
\ffl 1
\end{equation}
et le groupe d'inertie $I$ est lui-m\^{e}me objet d'un d\'{e}vissage
\begin{equation}\label{DevIn}
1\ffl P \ffl I \ffl I_{\rm mod}
\ffl 1
\end{equation}
o\`{u} $P$ est l'unique (pro-)$p$-sylow de $I$, et $I_{\rm mod}$ l'\og
inertie mod\'{e}r\'{e}e\fg. Ce dernier groupe s'identifie
canoniquement \`{a}
$\widehat{\ZZ}(1)(k_{\rm s})=\varprojlim_{n\geq1}\,\mu_{n}(k_{\rm
s})$; l'isomorphisme respecte les actions
de $\gal_{K}$ (par conjugaison sur $I_{\rm mod}$, par
l'action sur les racines de l'unit\'{e} sur $\widehat{\ZZ}(1)(k_{\rm
s})$). Le groupe $\widehat{\ZZ}(1)(k_{\rm s})$ s'identifie aussi
canoniquement \`{a} la partie
premi\`{e}re \`{a} $p$ de $\widehat{\ZZ}(1)(K_{\rm s})$.

On a donc une suite de sous-groupes de $\gal_{K}$, et la suite
correspondante de sous-corps de $K_{\rm s}$:
\begin{equation}\label{EqSsCorps}
\begin{array}{ccccc}
	\gal_{K} & \supset & I & \supset & P \cr
	K & \subset & K_{\rm nr} & \subset & K_{\rm mod}.
\end{array}
\end{equation}
Pour d\'{e}montrer le th\'{e}or\`{e}me \ref{ThRec}, nous appliquerons
\ref{ThPpal} en prenant pour suite exacte (\ref{DevPi}) la suite
\begin{equation}\label{DevPiMod}
1\ffl I_{\rm mod} \ffl \gal(K_{\rm mod}/K) \ffl  \gal_{k} \ffl 1.
\end{equation}

La preuve reposera sur les trois lemmes qui suivent, et qui montrent
respectivement que les trois conditions de \ref{ThPpal} sont
satisfaites.

\begin{lem}\label{ModScind} La suite exacte {\rm(\ref{DevPiMod})} de
groupes profinis est scind\'{e}e.
\end{lem}
\dem Soit $\varpi$ une uniformisante de $\Lambda$, et soit $L\subset
K_{\rm s}$ le corps de d\'{e}composition de tous les polyn\^{o}mes
$X^n-\varpi$, o\`{u} $n$ parcourt les entiers premiers \`{a} $p$.
On v\'{e}rifie alors sans mal que $L/K$ est totalement et
mod\'{e}r\'{e}ment ramifi\'{e}e, donc
$L\subset K_{\rm mod}$ et $L\cap K_{\rm nr}=K$. D'autre part il est
bien connu que $L\,K_{\rm nr}=K_{\rm mod}$: voir par exemple
\cite{C-F}, I, \S~8, Corollary 1 of
Proposition 1
(c'est l\`{a} un avatar du \og lemme d'Abhyankar\fg). Ceci
implique que le sous-groupe $\gal(K_{\rm mod}/L)$ de $\gal(K_{\rm
mod}/K)$ s'envoie isomorphiquement sur $\gal(K_{\rm nr}/K)=\gal_{k}$,
donc scinde la suite exacte de l'\'{e}nonc\'{e}.\qed

\begin{lem}\label{QuotImod} Soient $F$ un corps, $F_\mathrm{s}$ une
	cl\^{o}ture s\'{e}parable de $F$, et $A$ un $F$-sch\'{e}ma en
groupes fini \'{e}tale
	isomorphe (comme groupe avec action de $\gal(F_\mathrm{s}/F)$) \`{a}
un
	quotient de $\widehat{\ZZ}(1)(F_{\rm s})$. Alors,
$\ree(F,A)=\HH^1(F,A)$.
\end{lem}
\dem
L'hypoth\`{e}se implique que $A$ est isomorphe \`{a} $\mu_{n,F}$,
pour un entier $n$ premier \`{a} $\car(F)$. La th\'{e}orie de Kummer implique
imm\'{e}diatement que tout torseur sous un tel groupe est
\'{e}l\'{e}mentairement R-\'{e}quivalent au torseur trivial (c'est
d'ailleurs un cas par\-ti\-cu\-lier
de \ref{CritTVR}\,\ref{CritTVR2}).\qed
\medskip

\begin{lem}\label{RelevREq} Soient $\Lambda$, $K$ et $k$ comme en
\rref{NotVal}. Soit $G$ un $\Lambda$-sch\'{e}ma en groupes fini
\'{e}tale d'ordre inversible dans $\Lambda$, et soient
$\mathcal{X}$ et $\mathcal{Y}$ deux $G$-torseurs sur $\Spec\Lambda$.

Si les $G_k$-torseurs $\mathcal{X}_k$ et $\mathcal{Y}_k$ sont
\'{e}l\'{e}mentairement R-\'{e}quivalents (resp. R-\'{e}quivalents),
alors il en est de m\^{e}me des $G_K$-torseurs
$\mathcal{X}_K$ et $\mathcal{Y}_K$. \emph{(Bien entendu la notation
$G_{k}$ d\'{e}signe $G\times_{\Spec(\Lambda)}\Spec(k)$, etc.)}
\end{lem}
\dem Il suffit de traiter le cas de la R-\'{e}quivalence
\'{e}l\'{e}mentaire.
Il existe alors un rev\^{e}tement ramifi\'{e} $f_0:C_0\to\PP^1_k$,
o\`{u} $C_0$
est une $k$-courbe projective lisse munie d'une action de $G_k$ qui
fait de $f_0$ un $G_k$-torseur au-dessus d'un  ouvert de $\PP^1_k$
contenant $0$ et $1$, de telle sorte que
$f_0^{-1}(0)\cong\mathcal{X}_k$ et $f_0^{-1}(1)\cong\mathcal{Y}_k$
comme $G_k$-torseurs. Vu l'hypoth\`{e}se sur l'ordre de $G$, $f_0$
est
automatiquement mod\'{e}r\'{e}ment ramifi\'{e}; soit
$D_0\subset\PP^1_k$ son lieu
de  ramification, qui est un sous-sch\'{e}ma de $\PP^1_k$ \'{e}tale
sur
$\Spec k$. Choisissons un diviseur $D\subset\PP^1_\Lambda$ \'{e}tale
sur
$\Spec\Lambda$ et relevant $D_0$. La th\'{e}orie des d\'{e}formations
des
rev\^{e}tements mod\'{e}r\'{e}s (\cite{fulton}, th. 4.8, ou
\cite{revP1}, prop.
7.2.7) implique l'existence d'un unique rev\^{e}tement mod\'{e}r\'{e}
$f:C\to\PP^1_\Lambda$, relevant $f_0$, dont le diviseur de
ramification est $D$. Vu l'unicit\'{e}, $f$ est automatiquement un
$G$-rev\^{e}tement, et $f^{-1}(0)$ (resp. $f^{-1}(1)$) est un
$G$-torseur
sur $\Spec\Lambda$ relevant $\mathcal{X}_k$ (resp. $\mathcal{Y}_k$)
donc isomorphe \`{a} $\mathcal{X}$ (resp. \`{a} $\mathcal{Y}$), ce
qui ach\`{e}ve
la d\'{e}monstration.\qed

\Subsection{Preuve du th\'{e}or\`{e}me \ref{ThRec}. } Outre les
notations ci-dessus, on se donne
un groupe fini $G$ d'ordre premier \`{a} $p$. On suppose que
$\re(k,G)=\HH^1(k,G)$, et l'on veut en d\'{e}duire que
$\re(K,G)=\HH^1(K,G)$.

Comme on l'a annonc\'{e} plus haut, on applique \ref{ThPpal} en
prenant pour $M$ le corps $K_{\rm mod}$, et pour suite exacte
(\ref{DevPi}) la suite (\ref{DevPiMod}). Celle-ci  est bien
scind\'{e}e d'apr\`{e}s le lemme \ref{ModScind}. Noter que
\emph{tout} $G$-torseur sur $K$ est mod\'{e}r\'{e}ment ramifi\'{e},
donc trivialis\'{e} par $M$. L'extension $M_0$ de \ref{ThPpal} n'est
autre que $K_\mathrm{nr}$.

La condition \ref{ThPpal2} de \ref{ThPpal} est v\'{e}rifi\'{e}e
d'apr\`{e}s le
lemme \ref{QuotImod}. V\'{e}rifions la condition \ref{ThPpal1}: il
s'agit donc de
voir que tout $G$-torseur $X\to\Spec K$ \emph{non ramifi\'{e}} est
dans $\re(K,G)$. Un tel torseur se prolonge canoniquement en un
$G$-torseur $\mathcal{X}\to\Spec\Lambda$, de fibre sp\'{e}ciale un
$G$-torseur
${X_0}\to\Spec k$. Par hypoth\`{e}se, $X_0$ est dans $\re(k,G)$. Le
lemme \ref{RelevREq}
permet donc de conclure.\qed

\begin{rem} Le lemme \ref{RelevREq} n'appara\^{\i}t pas dans
\cite{gille}. Lorsque le corps r\'{e}siduel $k$ est fini, Gille
utilise
un autre argument pour v\'{e}rifier (au langage pr\`{e}s) la
condition
\ref{ThPpal1} de \ref{ThPpal}: puisque $\gal_k\cong\widehat{\ZZ}$,
l'\'{e}tude des $G$-torseurs non ramifi\'{e}s se ram\`{e}ne \`{a}
celle des
$\ZZ/n\ZZ$-torseurs, o\`{u} $n$ est premier \`{a} $p$, et un tel
torseur est
\'{e}l\'{e}mentairement R-\'{e}quivalent au torseur trivial d'apr\`{e}s
\ref{CritTVR}\,\ref{CritTVR2}.

C'est encore le fait que
$\gal_k\cong\widehat{\ZZ}$, plut\^{o}t que le lemme \ref{ModScind},
qu'utilise Gille pour voir que la suite (\ref{DevPi}) est scind\'{e}e.
\end{rem}

\begin{rem}
	L'exemple de Colliot-Th\'{e}l\`{e}ne cit\'{e} plus haut
(\ref{ContrExZ/8}) montre
	encore que l'hypoth\`{e}se sur l'ordre du groupe dans \ref{RelevREq}
	ne peut \^{e}tre supprim\'{e}e.
\end{rem}

\section{Cas d'un corps r\'{e}siduel ab\'{e}lien}\label{SecResAb}
\Subsection{Notations. }\label{NotSecResAb} On reprend les
hypoth\`{e}ses et notations de \ref{NotVal}. Si $A$ et $B$ sont deux
sous-groupes ferm\'{e}s de $\gal_{K}$, on
conviendra de noter $[A,B]$ le sous-groupe \emph{ferm\'{e}\/} de
$\gal_{K}$ engendr\'{e} par les commutateurs $[a,b]$ ($a\in A$, $b\in
B$). On pose alors
\begin{eqnarray}\label{DefP}
P_{0}&:=&P\cap[I,\gal_{K}]=\ker(I\varfl{{\rm can}}I/[I,\gal_{K}]\times I_{\rm
mod})\cr
M&:=&\vphantom{\frac{1}{2}}(K_{\rm s})^{P_{0}}.
\end{eqnarray}
Le diagramme (\ref{EqSsCorps}) se compl\`{e}te donc en:
$$
\begin{array}{ccccccccccc}
	\gal_{K} & \supset & I & \supset & P & \supset &
	P\cap[\gal_{K},\gal_{K}] & \supset & P_{0} & \supset &
	[I,I]\cr
	K & \subset & K_{\rm nr} & \subset & K_{\rm mod} & \subset & K_{\rm
	mod}K_{\rm ab} & \subset & M & \subset &
	(K_{\rm nr})_{\rm ab}\vphantom{{\displaystyle{\frac{1}{2}}}}
\end{array}
$$
o\`{u} l'on remarque que $[I,I]\subset P$ puisque $I_{\rm mod}$
est commutatif, d'o\`{u} l'inclusion $[I,I]\subset P_{0}$.

En passant au quotient par $P_{0}$ la suite (\ref{DevGal}), on
obtient une suite exacte de groupes profinis
\begin{equation}\label{DevPiBis}
1\ffl \Gamma \ffl \Pi \ffl \pi \ffl 1
\end{equation}
avec:
$$\begin{array}{rccll}
\Gamma & = & I/P_{0} & = & {\rm Im\,}(I\varfl{{\rm can}}I/[I,\gal_{K}]
\times I_{\rm mod})\cr
\Pi & = & \rlap{$\gal\,(M/K)$}\phantom{I/P_{0}}\cr
\pi & = & \gal_{k} & = & \gal\,(K_{\rm nr}/K).
\end{array}
$$

\begin{rem}\label{RemDevPi}
$\Gamma$ s'identifie \`{a} un sous-groupe de
$I/[I,\gal_{K}]\times\prod_{l\neq p}{\ZZ}_{l}(1)(K_{\rm s})$; il est
donc commutatif (ceci \'{e}quivaut d'ailleurs \`{a} la
propri\'{e}t\'{e} $[I,I]\subset
P_{0}$, d\'{e}j\`{a} observ\'{e}e). De plus, le noyau
$U:=\gal\,(K_{\rm s}/K_{\rm cycl})$
du caract\`{e}re cyclotomique op\`{e}re trivialement sur $\Gamma$;
autrement dit, l'image $U/P_{0}$ de $U$ dans $\Pi$ centralise
$\Gamma$. Ceci \'{e}quivaut encore \`{a} dire que $[I,U]\subset
P_{0}$, et l'on peut v\'{e}rifier
qu'en fait $P_{0}=[I,I][I,U]=[I,IU]$.
\end{rem}

\begin{thm}\label{ThResAb} Avec les notations de  \rref{NotSecResAb},
on se donne un groupe fini $G$, et l'on suppose que:
\begin{romlist}
\item\label{ThResAb1} $\pi=\gal_k$ est commutatif;
\item\label{ThResAb2} la suite {\rm(\ref{DevPiBis})} est scind\'{e}e;
\item\label{ThResAb3} la condition \cyc{K,G} de \rref{DefCondCyc} est 
satisfaite.
\end{romlist}

Alors, on a $\HH^1(M/K,G)\subset\re(K,G)$ \emph{(et donc
$\HH^1(M/K,G)\subset\ree(K,G)$ d'apr\`{e}s \cite{ReqSim}, cf. la
remarque \ref{RemThPpal5})}.
\end{thm}

\Subsubsection{Remarques. }\label{RemThResAb}

\noindent(1) Les conditions  \ref{ThResAb1} et  \ref{ThResAb2} du
th\'{e}or\`{e}me sont notamment v\'{e}rifi\'{e}es lorsque le groupe
$\pi=\gal_{k}$ est un quotient sans torsion de $\widehat{\ZZ}$, et en
particulier lorsque $k$ est fini ou s\'{e}parablement clos (et aussi
lorsque $k=\CC(\!(t)\!)$, mais voir (2) ci-dessous).\smallskip

\noindent(2) Lorsque $G$ est d'ordre premier \`{a} $p$, on n'obtient
ici rien de mieux que le th\'{e}or\`{e}me \ref{ThRec}: en effet on voit
tout de suite que \cyc{K,G} implique \cyc{k,G}, de sorte que, d'apr\`{e}s
\ref{CorGalProCyc}, l'hypoth\`{e}se de \ref{ThRec} est
v\'{e}rifi\'{e}e. En
dehors de cette situation, le cas particulier le plus simple de
\ref{ThResAb} est \ref{CorResAb} ci-dessous, que l'on peut d'ailleurs
facilement extraire de la d\'{e}monstration de \cite{gille}.
\smallskip

\begin{subcor}\label{CorResAb}
	On suppose que $p$ est un nombre premier impair et que $K$ est une
	ex\-ten\-sion finie de $\QQ_{p}$. Alors tout $G$-torseur
mod\'{e}r\'{e}ment
	ramifi\'{e} est dans $\re(K,G)$.
\end{subcor}
\dem Les conditions  \ref{ThResAb1} et  \ref{ThResAb2} de \ref{ThResAb}
sont v\'{e}rifi\'{e}es (remarque \ref{RemThResAb}\,(1)), et  \ref{ThResAb3}
l'est aussi puisque $p>2$ (remarque \ref{RemCPCyc}\,\ref{RemCPCyc3}).\qed

\Subsection{D\'{e}monstration de \ref{ThResAb}. }

Nous supposerons que $p>1$, sinon la remarque \ref{RemThResAb}\,(2)
s'applique.

On applique le th\'{e}or\`{e}me \ref{ThPpal} en prenant pour suite
(\ref{DevPi}) la suite (\ref{DevPiBis}), qui est scind\'{e}e par
hypoth\`{e}se. Le corps $M_0$ de \ref{ThPpal} est encore
$K_\mathrm{nr}$.

V\'{e}rifions la condition \ref{ThPpal1} de \ref{ThPpal}. Le
raisonnement
est le m\^{e}me que dans la preuve de \ref{CorGalProCyc}: un
$G$-torseur trivialis\'{e} par  $M_0$ correspond \`{a} un morphisme
continu $\pi\to G$. Comme $\pi$ est commutatif par
hypoth\`{e}se, un tel torseur est induit par un torseur sous un
sous-groupe commutatif $H$ de $G$. Comme \cyc{K,H} est
v\'{e}rifi\'{e}e, on peut
appliquer \ref{CritTVR}\,\ref{CritTVR5}.

V\'{e}rifions la condition \ref{ThPpal2} de \ref{ThPpal}. Soit donc
$H$
un quotient fini de $\Gamma$, dont le groupe sous-jacent est un
sous-groupe de $G$. Alors $H$ est commutatif (\ref{RemDevPi}) et l'on
peut donc supposer que c'est un $l$-groupe, pour $l$ premier.

Si $l=p$, alors la structure de $\Gamma$ montre que $H$ est un
quotient de $I/[I,\gal_{K}]$; or l'action de $\gal_K$ sur
$I/[I,\gal_{K}]$ est triviale de sorte que $H$ est constant. La
conclusion r\'{e}sulte donc de \ref{CritTVR}\,\ref{CritTVR5}.

Si $l\neq p$, alors $H$ est un sous-quotient de $A\times\ZZ_l(1)$,
o\`{u}
$A$ est la composante $l$-primaire de $I/[I,\gal_{K}]$ ($A$ est donc,
comme ci-dessus, un sch\'{e}ma en groupes constant). Donc $H$ est de
type
multiplicatif, d\'{e}ploy\'{e} par $K(\mu_{l^n})$ pour $n$ assez
grand. C'est
l\`{a} une extension \emph{cyclique} de $K$: c'est clair si $l\neq2$,
et
si $l=2$ alors $p>2$ de sorte que la remarque
\ref{RemCPCyc}\,\ref{RemCPCyc3} s'applique (on utilise ici
l'hypoth\`{e}se $p>1$ faite au d\'{e}but). On conclut par
\ref{CritTVR}\,\ref{CritTVR2}.\qed

\section{Bitorseurs: g\'{e}n\'{e}ralit\'{e}s}
\label{SecBitors}
On rappelle ci-dessous les principales propri\'{e}t\'{e}s des
bitorseurs; pour
plus de d\'{e}tails, voir par exemple \cite{giraud} ou \cite{breen}.
\Subsection{D\'{e}finition des bitorseurs.}
Soient $G'$ et  $G$ deux groupes d'un topos $T$ (par exemple deux
sch\'{e}mas en
groupes sur un corps). Rappelons (\cite{giraud}, \cite{breen})
qu'un $(G',G)$-\emph{bitorseur\/} est la donn\'{e}e d'un objet $X$ de
$T$,
muni d'une action \`{a} gauche de $G'$ et d'une action \`{a} droite
de $G$,
ces deux actions faisant respectivement de  $X$ un $G'$-torseur \`{a}
gauche et un $G$-torseur \`{a} droite, et de plus commutant entre
elles. En g\'{e}n\'{e}ral un tel bitorseur sera not\'{e} $(G',X,G)$,
sans
notation particuli\`{e}re pour les actions.

On voit imm\'{e}diatement, dans ces conditions, que $G'$ (resp.~$G$)
s'identifie
canoniquement au faisceau $\soul{\rm Aut}_{G}(X)$ des
$G$-automorphismes de $X$
(resp. \`{a} $\soul{\rm Aut}_{G'}(X)^{\circ}$, groupe oppos\'{e} au
faisceau des $G'$-automorphismes de $X$), de sorte que, par exemple,
tout
$G$-torseur \`{a} droite $Y$ d\'{e}termine canoniquement un bitorseur
$(\soul{\rm Aut}_{G}(Y),Y, G)$,
et que tout bitorseur est de cette forme, \`{a} isomorphisme
canonique pr\`{e}s (voir ci-dessous pour la notion d' isomorphisme).
\Subsection{Morphismes de bitorseurs. }\label{MorBitors}Soient
${\cal X}_{1}=(G'_{1},X_{1},G_{1})$ et ${\cal
X}_{2}=(G'_{2},X_{2},G_{2})$ deux
bitorseurs. Nous appellerons \emph{morphisme} de ${\cal X}_{1}$ vers
${\cal X}_{2}$ un triplet $\Phi=(\varphi':G'_{1}\to G'_{2},
u:X_{1}\to
X_{2}, \varphi:G_{1}\to G_{2})$ de morphismes de $T$, tel que
$\varphi'$ et $\varphi$ soient des morphismes de groupes et que $u$
soit $\varphi'$-\'{e}quivariant \`{a} gauche et
$\varphi$-\'{e}quivariant \`{a} droite.

On obtient ainsi une cat\'{e}gorie not\'{e}e ${\rm Bitors}\,(T)$.

Il arrive que l'on ait besoin d'imposer les morphismes
$\varphi',\varphi$ (avec les notations ci-dessus), ou
seulement l'un d'eux: on dira donc que $\Phi$ est un
$(\varphi',\varphi)$-morphisme, ou un
$(\varphi',\ast)$-morphisme,
ou un $(\ast,\varphi)$-morphisme. De m\^{e}me, on parlera de
$(G',\ast)$-bitorseurs, et de $(\ast,G)$-bitorseurs.

Si $\Phi=(\varphi', u,\varphi)$ est un morphisme de bitorseurs, il
est \'{e}quivalent de dire que $\varphi'$ est injectif (resp.
surjectif)
comme morphisme de $T$, ou que $u$ l'est, ou que $\varphi$ l'est;
nous dirons alors que $\Phi$ est injectif (resp. surjectif). De
m\^eme, $\Phi$ est un isomorphisme dans ${\rm Bitors}\,(T)$ si et
seulement si $\varphi'$ (ou $u$, ou $\varphi$) en est un dans $T$.

Tout morphisme $\Phi$ s'\'{e}crit de mani\`{e}re essentiellement unique
sous la forme $\beta\circ\alpha$, o\`{u} $\beta$ est injectif et $\alpha$
surjectif. Bien entendu, le bitorseur source de $\beta$ n'est autre,
\`{a} isomorphisme canonique pr\`{e}s, que
$({\rm Im}\,\varphi', {\rm Im}\,u, {\rm Im}\,\varphi)$; il
sera appel\'{e} l'\emph{image} de $\Phi$.
\Subsection{Trivialisations. }\label{SubSecTriv}
Si $G$ est un groupe de $T$, le \emph{$G$-bitorseur trivial\/} est
par
d\'{e}finition le $(G,G)$-bitorseur ${\rm Triv}(G)=(G,G,G)$ o\`{u}
$G$
op\`{e}re \`{a} droite et \`{a} gauche sur lui-m\^{e}me par
translations.

Si $(G',X,G)$ est un bitorseur, il revient au m\^{e}me de dire que le
$G'$-torseur \`{a} gauche $X$ est trivial, ou que le
$G$-torseur \`{a} droite $X$ est trivial, ou que l'objet $X$ de $T$ a
une
section; de plus ces conditions sont v\'{e}rifi\'{e}es localement
dans $T$.

Si ces conditions sont v\'{e}rifi\'{e}es, nous commettrons l'abus
(courant
d\'{e}j\`{a} pour les torseurs) de dire que \og le bitorseur
$(G',X,G)$ est
trivial\fg.

De fa\c{c}on plus pr\'{e}cise, le choix d'une section $x$ de $X$
d\'{e}termine un
isomorphisme de $G$-torseurs \`{a} droite $u:G\to X$, donn\'{e} par
$g\mapsto xg$. Celui-ci se
prolonge en un $(\ast,{\rm Id}_{G})$-isomorphisme de la forme $({\rm
conj}\,(x),u,{\rm Id}_{G}): {\rm Triv}(G) \to (G',X,G)$, o\`{u}
${\rm conj}\,(x):G\to G'$ est un isomorphisme de groupes de $T$ qui
m\'{e}rite le
nom de \og conjugaison par $x$\fg; de fait, si $(G',X,G)$ est le
bitorseur trivial $(G,G,G)$, on v\'{e}rifie tout de suite que
${\rm conj}\,(x):G\to G$ est l'automorphisme int\'{e}rieur $g\mapsto
xgx^{-1}$.

\Subsection{Bitorseurs et classes de conjugaison. }\label{ClasConj}
La propri\'{e}t\'{e} qui pr\'{e}c\`{e}de (que l'on exprime couramment
en
disant que $G'$ est une \og forme int\'{e}rieure\fg\ de $G$) a une
cons\'{e}quence importante: \`{a} tout bitorseur $(G',X,G)$ est
associ\'{e}e une
bijection canonique entre l'ensemble
des sous-objets de $G$ invariants par conjugaison et l'ensemble
des sous-objets de $G'$ invariants par conjugaison. En particulier,
si $H$ est un \emph{sous-groupe distingu\'{e}\/} de $G$, il lui
correspond
canoniquement un sous-groupe distingu\'{e} $H'$ de $G'$,
caract\'{e}ris\'{e}
par la propri\'{e}t\'{e} que $X/H=H'\backslash X$; noter que l'on a
alors un
morphisme canonique $(G',X,G)\to(G'/H',X/H,G/H)$.
\Subsection{Changement de groupe structural.
}\label{ChangGroupe}Soient
$\mathcal{X}=(G'_{1},X_{1},G_{1})$ un bitorseur et $\varphi:G_{1}\to
G_{2}$ un
morphisme de groupes de $T$. Il existe alors un $(\ast,G_{2})$-bitorseur
${\cal X}^{\varphi}=(G'_{2},X_{2},G_{2})$ et un
$(\ast,\varphi)$-morphisme $\Phi=(\varphi',u,\varphi):{\cal
X}\to{\cal X}^{\varphi}$ qui est universel au sens suivant: tout
$(\ast,\varphi)$-morphisme
${\cal X}\to{\cal Y}$ se factorise de fa\c{c}on unique sous la
forme    ${\cal X}\sfl{\Phi}{\cal X}^{\varphi}\sfl{\Theta}{\cal Y}$
o\`{u} $\Theta$ est un $(\ast,{\rm Id}_{G_{2}})$-(iso)morphisme.

(En d'autres termes, le foncteur $(G',X,G)\mapsto G$ fait de ${\rm
Bitors\,}(T)$ une cat\'{e}gorie cofibr\'{e}e en groupo\"{\i}des
au-dessus
de la cat\'{e}gorie des groupes de $T$.)

En tant que $G_{2}$-torseur \`{a} droite, $X_{2}$ n'est autre que le
torseur d\'{e}duit de $X_{1}$ par le changement de groupe structural
$\varphi$: c'est donc le produit contract\'{e} $X_{1}\times^{G_{1}}G_{2}$,
quotient de $X_{1}\times G_{2}$ par l'action de $G_{1}$ donn\'{e}e par
$((x_{1},g_{2}),g_{1})\mapsto(x_{1}g_{1},\varphi(g_{1})^{-1}g_{2})$.

Noter que dans cette construction, le groupe $G'_{2}$ et le
morphisme $\varphi':G'_{1}\to G'_{2}$ d\'{e}\-pen\-dent de $X_{1}$
(et pas
seulement de $\varphi$ et de $G'_{1}$). Cependant, $\varphi'$ est \og
localement isomorphe \`{a} $\varphi$\fg: plus pr\'{e}cis\'{e}ment, le
choix d'une
trivialisation de $\mathcal{X}$ d\'{e}termine des isomorphismes
$\psi_{1}:G_{1}\flis G'_{1}$ et $\psi_{2}:G_{2}\flis G'_{2}$ tels
que $\varphi'\psi_{1}=\psi_{2}\varphi$. On voit ainsi, notamment, que
si
$(G',X,G)$ est un bitorseur et $H$ un
sous-groupe de $G$ tel que $X$ soit induit (comme $G$-torseur \`{a}
droite) par un $H$-torseur, alors il est induit comme $G'$-torseur
\`{a}
gauche
par un torseur sous un sous-groupe $H'$ de $G'$, localement isomorphe
\`{a} $H$.

Enfin, on a bien entendu une op\'{e}ration similaire du c\^{o}t\'{e}
gauche: si
$\varphi':G'_{1}\to G'_{2}$ est un morphisme de groupes, elle fournit
un $(G'_{2},\ast)$-bitorseur not\'{e} ${}^{\varphi'}\!{\cal X}$ et un
$(\varphi',\ast)$-morphisme universel ${\cal X}\to{}^{\varphi'}\!{\cal X}$.
\Subsection{Cas particulier: passage au quotient.
}\label{CPQuotient}Soient
${\cal X}=(G',X,G)$ un bitorseur, $H\sinj{j}G$ un sous-groupe
\emph{distingu\'{e}\/} de
$G$, et $\varphi:G\to G/H$ le morphisme canonique. On a alors,
d'apr\`{e}s
\ref{SubSecTriv}, un sous-groupe distingu\'{e} $H'$ de $G'$ et un
morphisme ${\cal X}\to(G'/H',X/H,G/H)$. Bien entendu, celui-ci n'est
autre que le morphisme ${\cal X}\to{\cal X}^\varphi$ d\'{e}fini
ci-dessus.
Le bitorseur $(G'/H',X/H,G/H)$ pourra \^{e}tre not\'{e},
indiff\'{e}remment,
${\cal X}/H$ ou $H'\backslash {\cal X}$. On v\'{e}rifie, sans
surprise, la
propri\'{e}t\'{e} suivante:

\begin{sublem}\label{LemQotient} Avec les notations de
\rref{CPQuotient}, les conditions suivantes sont \'{e}quivalentes:
	\begin{romlist}
		\item le bitorseur ${\cal X}/H$ est trivial;
		\item le $G$-torseur \`{a} droite $X$ est induit par 
un $H$-torseur;
		\item le $G'$-torseur \`{a} gauche $X$ est induit par un
$H'$-torseur;
		\item il existe un $(H',H)$-bitorseur $Y$ et un
$(\ast,j)$-morphisme
		$(H',Y,H)\to {\cal X}$.\qed
	\end{romlist}
\end{sublem}
\Subsection{Composition de bitorseurs. }La cat\'{e}gorie ${\rm
Bitors\,}(T)$
poss\`{e}de un \og produit\fg\ partiellement d\'{e}fini,
qui fait tout son charme: si ${\cal X}_{1}=(G_{1},X_{1},G_{2})$ et
${\cal X}_{2}=(G_{2},X_{2},G_{3})$ sont dans ${\rm Bitors\,}(T)$, le
produit
contract\'{e} $X_{1}\times^{G_{2}}X_{2}$ admet une structure
naturelle
de $(G_{1},G_{3})$-bitorseur. On notera
$$
{\cal X}_{1}\wedge{\cal X}_{2}=
(G_{1},X_{1}\times^{G_{2}}X_{2},G_{3})
$$
le bitorseur ainsi obtenu. La loi $\wedge$ admet une contrainte
d'associativit\'{e} (des isomorphismes $({\cal X}_{1}\wedge{\cal
X}_{2})
\wedge{\cal X}_{3}\flis{\cal X}_{1}\wedge({\cal X}_{2}
\wedge{\cal X}_{3})$ assortis de compatibilit\'{e}s); les bitorseurs
triviaux sont neutres \`{a} droite et \`{a} gauche, et tout
bitorseur  ${\cal
X}=(G',X,G)$ admet un inverse, \`{a} savoir  ${\cal
X}^{-1}=(G,X^{-1},G')$, o\`{u}
$X^{-1}=X$ muni de l'action \`{a} gauche (resp. \`{a} droite) de $G$
(resp. $G'$)
donn\'{e}e par $(g,x)\mapsto xg^{-1}$ (resp. $(x,g')\mapsto
{g'}^{-1}x$).

\begin{sublem}\label{MorProd}
	Soient  ${\cal X}_{1}=(G_{1},X_{1},G_{2})$,
${\cal X}_{2}=(G_{2},X_{2},G_{3})$, $\mathcal{Y}=(H',Y,H)$ des
bi\-tor\-seurs (avec ${\cal X}_{1}$ et ${\cal X}_{2}$ composables),
et
$$\Phi:{\cal X}_{1}\wedge{\cal X}_{2}\ffl \mathcal{Y}$$
un morphisme.

Il existe alors un groupe $G'_{2}$ de $T$, et un morphisme de groupes
$\varphi_{2}:G_{2}\to G'_{2}$ tels que $\Phi$ se factorise sous la
forme
$${\cal X}_{1}\wedge{\cal X}_{2}\varfl{\Phi_{1}\wedge\Phi_{2}}{\cal
X}_{1}^{\varphi_{2}}\wedge{}_{\vphantom{1}}^{\varphi_{2}}\!{\cal
X}_{2}\varfl{\Psi}\mathcal{Y}
$$
o\`{u} $\Phi_{1}:{\cal X}_{1}\to{\cal X}_{1}^{\varphi_{2}}$ et
$\Phi_{2}:{\cal X}_{2}\to{}_{\vphantom{1}}^{\varphi_{2}}\!{\cal 
X}_{2}$ sont les
morphismes canoniques, et o\`{u} $\Psi$ est un isomorphisme.
\end{sublem}
\dem \'Ecrivons $\Phi=(\varphi_{1},u,\varphi_{3}):{\cal
X}_{1}\wedge{\cal
X}_{2}\fl(H',Y,H)$. Consid\'{e}rons le morphisme canonique
$\Phi_{2}:{\cal
X}_{2}\to{\cal X}_{2}^{\varphi_{3}}$: il est de la forme
$(\varphi_{2},u_{2},\varphi_{3})$ o\`{u} $\varphi_{2}:G_{2}\to
G'_{2}$
est un morphisme de groupes; le bitorseur ${\cal
X}_{2}^{\varphi_{3}}$ s'identifie aussi \`{a} 
${}_{\vphantom{1}}^{\varphi_{2}}\!{\cal
X}_{2}$.
Le morphisme canonique ${\Phi_{1}\wedge\Phi_{2}}:{\cal
X}_{1}\wedge{\cal
X}_{2}\to {\cal 
X}_{1}^{\varphi_{2}}\wedge{}_{\vphantom{1}}^{\varphi_{2}}\!{\cal
X}_{2}$ est par
construction un $(\ast,\varphi_{3})$-morphisme, comme $\Phi$, de
sorte
que $\mathcal{Y}$ et
${\cal X}_{1}^{\varphi_{2}}\wedge{}_{\vphantom{1}}^{\varphi_{2}}\!{\cal X}_{2}$
s'identifient tous deux \`{a} $({\cal X}_{1}\wedge{\cal
X}_{2})^{\varphi_{3}}$. Le lemme en r\'{e}sulte.\qed

\Subsection{Bitorseurs d'isomorphismes. }Soit $G$ un groupe de $T$;
consid\'{e}rons deux $(\ast, G)$-bitorseurs $(G',X,G)$ et
$(G'',Y,G)$. Alors
$\soul{\rm Isom}_{G}(X,Y)$ admet une structure naturelle de
$(G'',G')$-bitorseur,
et l'on a un isomorphisme canonique de $(G'',G')$-bitorseurs
\begin{equation}
	\label{EqIsom}
	(G'',\soul{\rm Isom}_{G}(X,Y),G')\;\flis\;
(G'',Y,G)\wedge(G',X,G)^{-1}.
\end{equation}
En combinant ceci avec le lemme \ref{LemQotient}, on obtient le
r\'{e}sultat
suivant, qui nous servira plus loin:

\begin{sublem}\label{LemIsomQuot}
	Soient ${\cal X}=(G',X,G)$, ${\cal Y}=(G'',Y,G)$ deux
	$(\ast, G)$-bitorseurs, $H$ un sous-groupe distingu\'{e} de $G$,
	$H'\subset G'$ et $H''\subset G''$ les sous-groupes
	qui lui  correspondent
	ca\-no\-ni\-que\-ment
	via ${\cal X}$ et ${\cal Y}$ respectivement (cf. \rref{SubSecTriv}).
	On suppose que ${\cal X}/H$ et ${\cal Y}/H$ sont isomorphes. Alors
	le $(G'', G')$-bitorseur ${\cal Y}\wedge{\cal X}^{-1}$ provient par
	changement de groupe structural d'un $(H'',H')$-bitorseur.\qed
\end{sublem}

\section{Bitorseurs et R-\'{e}quivalence}
\label{SecBitorsREquiv}
\Subsection{Notations. }\label{NotSecBitorsREquiv}Soit $K$ un corps,
dont on
fixe une cl\^{o}ture
s\'{e}parable $K_{\rm s}$. Nous allons appliquer les
notions ci-dessus au topos $K_{\rm et}$ des faisceaux sur le petit
site \'{e}tale de $K$ (qui est \'{e}quivalent au topos des
$\gals{K}$-ensembles).

Si $G$ est un $K$-sch\'{e}ma en groupes fini \'{e}tale, nous noterons
\begin{equation*}
	(X,G)\rree(Y,G),\quad\hbox{ resp.}\quad(G,X')\rree(G,Y')
\end{equation*}
la R-\'{e}quivalence \'{e}l\'{e}mentaire entre deux $G$-torseurs
\`{a} droite $X$ et
$Y$ (resp.~deux $G$-torseurs \`{a} gauche $X'$ et $Y'$). De m\^{e}me
nous
noterons
\begin{equation*}
	(X,G) \rre(Y,G),\quad\hbox{ resp.}\quad(G,X')\rre(G,Y')
\end{equation*}
la R-\'{e}quivalence. La notation $(G,G)$ d\'{e}signera le
$G$-torseur
trivial (\`{a} droite ou \`{a} gauche, suivant le contexte).

Dans ce qui suit, tous les sch\'{e}mas en groupes et (bi)torseurs
consid\'{e}r\'{e}s sont finis \'{e}tales sur $K$.

\begin{prop}
	\label{REquivElBit} {\rm(Bitorseurs et R-\'{e}quivalence
	\'{e}l\'{e}mentaire)}\par
	\begin{romlist}
		\item\label{REquivElBit1} Soient $H$, $G$, $G'$, 
$G''$, $G'''$ des
		$K$-sch\'{e}mas en groupes
		finis \'{e}tales, $(G',X,G)$, $(G'',Y,G)$, $(G,Z,G''')$ des
		bitorseurs et $\varphi:G\to H$ un morphisme de $K$-groupes.
		Alors on a les implications suivantes:
		\begin{subromlist}
			\item\label{REquivElBit11}
 
	$(X,G)\rree(Y,G)\;\Longleftrightarrow\;(G,X^{-1})\rree(G,Y^{-1})$;
			\item\label{REquivElBit12}
			$(X,G)\rree(Y,G)\;\Longrightarrow\;(X\times^G 
H, H)\rree(Y\times^G
H,H)$;
			\item\label{REquivElBit13}
			$(X,G)\rree(Y,G)\;\Longrightarrow\;
			(X\times^G Z, G''')\rree(Y\times^G Z, G''')$;
			\item\label{REquivElBit14} 
$(X,G)\rree(G,G)\;\Longleftrightarrow\;
			(G',X)\rree(G',G')$;
			\item\label{REquivElBit15}
			$\begin{array}[t]{rcl}
			(X,G)\rree(Y,G)
			& \Longleftrightarrow & (G'',\soul{\rm
Isom}_{G}(X,Y))\rree(G'',G'')\cr
			& \Longleftrightarrow & (\soul{\rm 
Isom}_{G}(X,Y),G')\rree(G',G')
		\end{array}
		$
	\end{subromlist}
	ainsi que l'analogue de \rref{REquivElBit13} pour le produit
	contract\'{e} \`{a}
	gauche, et les analogues de \rref{REquivElBit12} et
	\rref{REquivElBit15} pour les
	torseurs \`{a} gauche.
	\smallskip

	\item\label{REquivElBit2} Notons $\gree$ la classe des
bitorseurs $(G',X,G)$
	tels que $(X,G)\rree(G,G)$ (ou tels que
	$(G',X)\rree(G',G')$, ce qui revient au m\^{e}me d'apr\`{e}s
	\ref{REquivElBit1}\,\ref{REquivElBit14}). Alors $\gree$
poss\`{e}de les
	propri\'{e}t\'{e}s suivantes:
	\begin{subromlist}
		\item\label{REquivElBit21} $\gree$ contient les bitorseurs
triviaux
		${\rm Triv}(G)$;
		\item\label{REquivElBit22} $\gree$ est stable par inverse:
si ${\cal X}\in\gree$, alors ${\cal X}^{-1}\in\gree$;
		\item\label{REquivElBit23} $\gree$ est \og stable par
		morphismes\fg: si ${\cal X}\in\gree$ et si
		$\Phi:{\cal X}\to{\cal Y}$ est un morphisme, alors ${\cal
Y}\in\gree$.
	\end{subromlist}
	\smallskip

	\item\label{REquivElBit3} La R-\'{e}quivalence \'{e}l\'{e}mentaire
de torseurs est
	d\'{e}termin\'{e}e par la classe $\gree$ de
\rref{REquivElBit2}:
	de fa\c{c}on pr\'{e}cise (pour les torseurs
	\`{a} droite) si ${\cal X}=(G',X,G)$ et ${\cal Y}=(G'',Y,G)$ sont
deux
	$(\ast, G)$-bitorseurs,
	les propri\'{e}t\'{e}s suivantes sont \'{e}quivalentes:
	\begin{subromlist}
		\item\label{REquivElBit31} $(X,G)\rree(Y,G)$;
		\item\label{REquivElBit32} ${\cal Y}\wedge{\cal 
X}^{-1}\in\gree$;
		\item\label{REquivElBit33} ${\cal Y}$ est de la forme ${\cal
Z}\wedge{\cal X}$,
		avec ${\cal Z}\in\gree$.
	\end{subromlist}
\end{romlist}
\end{prop}
\dem \ref{REquivElBit1} Les assertions \ref{REquivElBit11},
\ref{REquivElBit12} et \ref{REquivElBit13} r\'{e}sultent
imm\'{e}diatement des
d\'{e}finitions.
Pour montrer \ref{REquivElBit14}, on applique \ref{REquivElBit13} en
prenant
$(G'',Y,G)=(G,G,G)$ et $(G,Z,G''')=(G,X^{-1},G')$: on obtient
$(G',G')\rree(X^{-1},G')$ d'o\`{u} $(G',X)\rree(G',G')$ d'apr\`{e}s
\ref{REquivElBit11}.

On en d\'{e}duit \ref{REquivElBit15} en utilisant (\ref{EqIsom}).

Finalement, \ref{REquivElBit2} (resp. \ref{REquivElBit3}) n'est
qu'une
reformulation des propri\'{e}t\'{e}s
\ref{REquivElBit11}, \ref{REquivElBit12} et \ref{REquivElBit13}
(resp. \ref{REquivElBit15}) de \ref{REquivElBit1}.\qed

\begin{prop}\label{PropREqBit} {\rm(Bitorseurs et
R-\'{e}quivalence)\/} Notons encore $\gree$ la classe de
bitorseurs d\'{e}finie en \rref{REquivElBit}\,\ref{REquivElBit2}.
	\begin{romlist}
		\item\label{PropREqBit1} Soit ${\cal X}=(G',X,G)$ un bitorseur.
		Les conditions suivantes sont \'{e}quivalentes:
		\begin{subromlist}
			\item\label{PropREqBit11} $(X,G)\rre(G,G)$;
			\item\label{PropREqBit12} $(G',X)\rre(G',G')$;
			\item\label{PropREqBit13} ${\cal X}$ est 
(\`{a} isomorphisme
			pr\`{e}s) de la forme
			${\cal X}_{1}\wedge{\cal 
X}_{2}\wedge\cdots\wedge{\cal X}_{n}$,
			o\`{u} les ${\cal X}_{i}$ sont dans $\gree$.
		\end{subromlist}
		\item\label{PropREqBit2} Soit $\gre$ la classe des bitorseurs
		${\cal X}=(G',X,G)$
		v\'{e}rifiant les conditions de \rref{PropREqBit1}. 
Alors $\gre$
v\'{e}rifie 		les analogues des propri\'{e}t\'{e}s de
		\rref{REquivElBit}\,\ref{REquivElBit2} (elle contient les
bitorseurs
		triviaux, et  est stable par inverse et par 
morphismes), et est de
		plus stable par la composition de bitorseurs (lorsqu'elle est
		d\'{e}finie).
		Plus pr\'{e}cis\'{e}ment, $\gre$ est la plus petite classe de
		bitorseurs stable par isomorphisme et par 
composition, et contenant
		$\gree$.
		\item\label{PropREqBit3} Soient ${\cal X}=(G',X,G)$ et ${\cal
		Y}=(G'',Y,G)$ deux $(\ast, G)$-bitorseurs. Les conditions
		sui\-vantes sont \'{e}quivalentes (o\`{u} $\gre$ est 
la classe de
		bitorseurs d\'{e}finie en \ref{PropREqBit2}):
		\begin{subromlist}
			\item\label{PropREqBit21} $(X,G)\rre(Y,G)$;
			\item\label{PropREqBit22} ${\cal 
Y}\wedge{\cal X}^{-1}\in{\gre}$;
			\item\label{PropREqBit23} ${\cal Y}$ est de 
la forme ${\cal
			Z}\wedge{\cal X}$, avec ${\cal Z}\in{\gre}$.
		\end{subromlist}
	\end{romlist}
\end{prop}
\dem L'assertion \ref{PropREqBit1} r\'{e}sulte de la d\'{e}finition
et de
\ref{REquivElBit}\,\ref{REquivElBit3}, et entra\^{\i}ne facilement
les autres.\qed

\begin{rem}\label{RemOrdreREquiv} On  peut \'{e}videmment
	pr\'{e}ciser \ref{PropREqBit1} de la fa\c{c}on suivante: si, pour
	$n\in\NN$
	donn\'{e}, l'on d\'{e}finit $\re_{n}(K,G)$ commme en
	\ref{RemThPpal3}, on a
	l'\'{e}quivalence: $(X,G)\in\re_{n}(K,G)\Leftrightarrow\mathcal{X}$
	est de la forme ${\cal X}_{1}\wedge{\cal
	X}_{2}\wedge\cdots\wedge{\cal
	X}_{n}$, o\`{u} les ${\cal X}_{i}$ sont dans $\gree$.
\end{rem}

\section{Bitorseurs \`{a} groupe structural constant}\label{SecBitorsGal}
Dans ce paragraphe nous d\'{e}crivons les bitorseurs dont l'un des deux
groupes structuraux est fini constant, lorsque le topos $T$ est galoisien,
c'est-\`{a}-dire \'{e}quivalent \`{a} la cat\'{e}gorie des ensembles avec
action continue d'un groupe profini.

\Subsection{$\Pi$-ensembles: notations et conventions.
}\label{NotPiEns}
On se donne d\'{e}sormais un groupe profini $\Pi$, et l'on
travaille dans la cat\'{e}gorie ${\cal C}_{\Pi}$ des $\Pi$-ensembles,
c'est-\`{a}-dire des ensembles $X$ munis d'une action \`{a} gauche continue de
$\Pi$ (pour la topologie profinie sur $\Pi$ et la topologie discr\`{e}te
sur $X$).
Cette cat\'{e}gorie est un topos (\cite{sga4}, IV, 2.4).

Si $\sf X$ est un $\Pi$-ensemble, son ensemble sous-jacent sera
not\'{e} $\vert{\sf X}\vert$;
le transform\'{e} de $x\in \vert {\sf X}\vert$ par $\sigma\in\Pi$
sera not\'{e}
$^{\sigma}\!x$.

Pour \'{e}viter des confusions, nous utiliserons (sauf exceptions
telles que la notation $\vert{\sf X}\vert$ ci-dessus) des
typographies diff\'{e}rentes pour les
ensembles ($X,Y,\ldots$) et les $\Pi$-ensembles ($\sf X,Y,\ldots$).

Si $\sf X$ et $\sf Y$ sont deux $\Pi$-ensembles \emph{finis}, l'objet
$\soul{\sf Hom}\,({\sf X},{\sf Y})$ de ${\cal C}_{\Pi}$ est
l'ensemble $H$ des
applications de $\vert {\sf X}\vert$ dans $\vert {\sf Y}\vert$ muni
de l'action
de $\Pi$ donn\'{e}e par la formule
$$(^\sigma\! f)(x)=^\sigma\![f(^{\sigma^{-1}}\!x)]$$
(pour $f\in H$, $\sigma\in\Pi$ et $x\in \vert {\sf X}\vert$).
La finitude assure que cette action est bien continue. (Sans
l'hypoth\`{e}se de finitude, la bonne d\'{e}finition de l'ensemble
sous-jacent \`{a} $\soul{\sf Hom}\,({\sf X},{\sf Y})$ est
$\varinjlim_{U}\,{\rm Hom}_{U}(|{\sf X}|,|{\sf Y}|)$ o\`{u} $U$
parcourt les sous-groupes ouverts de $\Pi$).

Les objets de ${\cal C}_{\Pi}$ pour lesquels l'action de $\Pi$ est
triviale seront dits \emph{constants\/} (ce sont les objets
constants du  topos ${\cal C}_{\Pi}$). L'objet constant associ\'{e}
\`{a} un
ensemble $X$ sera not\'{e} $\ssoul{X}$.

Un \emph{$\Pi$-groupe\/} est un groupe du topos ${\cal C}_{\Pi}$,
c'est-\`{a}-dire
un groupe muni d'une action \`{a} gauche
de $\Pi$ par automorphismes. Si ${\sf G}$ est un $\Pi$-groupe, on
notera encore $\vert {\sf G}\vert$, par abus, le \emph{groupe\/}
sous-jacent
\`{a} $\sf G$. Un $\sf G$-torseur \`{a} droite est alors un
$\Pi$-ensemble non vide
$\sf X$ muni d'une action \`{a} droite libre et transitive de $\vert
{\sf G}\vert$
sur $\vert{\sf X}\vert$ (not\'{e}e $(x,g)\mapsto xg$), la
compatibilit\'{e} avec les actions de $\Pi$ \'{e}tant donn\'{e}e par
la formule
$^\sigma\!(xg)=(^\sigma\! x)(^\sigma\! g)$ ($\sigma\in\Pi$, $x\in\vert{\sf
X}\vert$,
$g\in\vert{\sf G}\vert$).

On peut voir $\Pi$ lui-m\^{e}me (ainsi que ses sous-groupes
ferm\'{e}s
distingu\'{e}s et ses quotients) comme un pro-$\Pi$-groupe, en le
munissant de son
action par automorphismes int\'{e}rieurs.

Si ${\cal X}=({\sf G', X,G})$ est un bitorseur dans ${\cal C}_{\Pi}$,
on notera $|{\cal X}|$ le bitorseur ensembliste $({\sf |G'|,
|X|,|G|})$;
le foncteur ${\cal X}\mapsto|{\cal X}|$
d'oubli des actions de $\Pi$ est fid\`{e}le et compatible (notamment)
aux
changements de groupe structural et \`{a} la composition de
bitorseurs.

\Subsection{$\Pi$-ensembles: (bi)torseurs sous les groupes constants.
}
\label{SorConst}
\Subsubsection{Notations. }\label{NotConst}Consid\'{e}rons la
sous-cat\'{e}gorie
pleine ${\rm B}_{\Pi}$ de ${\rm Bitors\,}({\cal C}_{\Pi})$ form\'{e}
des
bi\-tor\-seurs
${\cal X}=({\sf G', X,G})$ tels que le $\Pi$-groupe ${\sf G}$ soit
\emph{constant}.

D\'{e}signons d'autre part par  ${\rm B}'_{\Pi}$ la cat\'{e}gorie
suivante:
\begin{itemize}
	\item un objet de  ${\rm B}'_{\Pi}$ est de la forme
	$(G',X,G,\theta)$, o\`{u} $(G',X,G)$ est un bitorseur ensembliste et
	o\`{u} $\theta:\Pi\to G'$ est un homomorphisme continu;
	\item un morphisme de $(G',X,G,\theta)$ vers $(H',Y,H,\psi)$ est un
	morphisme de bitorseurs $(\varphi',u,\varphi):(G',X,G)\to(H',Y,H)$
	v\'{e}rifiant $\psi\circ\varphi'=\theta:\Pi\to H'$.
\end{itemize}

\Subsubsection{Remarques. }\label{RemConst1}
La cat\'{e}gorie  ${\rm B}_{\Pi}$  est \'{e}quivalente \`{a}
la cat\'{e}gorie, en apparence plus simple, des couples $({\sf X},G)$
o\`{u}
$G$ est un groupe et ${\sf X}$
un torseur \`{a} droite sous le groupe constant $\ssoul{G}$:
l'\'{e}quivalence est
donn\'{e}e par le foncteur $({\sf G', X,G})\mapsto({\sf X,|G|})$, de
quasi-inverse  $({\sf X},G)\mapsto({\sf Aut}_{\ssoul{{\scriptstyle G}}}({\sf
X}),{\sf
X},\ssoul{G})$.

De la m\^{e}me fa\c{c}on, ${\rm B}'_{\Pi}$ est \'{e}quivalente \`{a}
la cat\'{e}gorie des objets $(X,G',\theta)$ o\`{u}
$G'$ est un groupe, $X$ un $G'$-torseur \`{a} droite, et
$\theta:\Pi\to G'$ un
homomorphisme continu.

Les d\'{e}finitions en termes de bitorseurs expriment mieux la
sym\'{e}trie de
la situation, et notamment le fait que l'on a un diagramme de \og
foncteurs d'oubli\fg:

$$\begin{array}{cccc}
& & {\rm B}'_{\Pi} & (G',X,G,\theta) \cr
& & \mapdown{\omega'} & \longmapstodown \cr
{\rm B}_{\Pi} & \sffl{\omega} & \underset{\mathstrut}{\rm 
Bitors\,(Ens)} & (G',X,G) \cr
{\cal X} & \longmapsto & |{\cal X}| \cr
\end{array}$$

\begin{subprop}\label{EquivBitorsConst}
	Il existe une \'{e}quivalence de cat\'{e}gories
	$$\Phi:{\rm B}_{\Pi}\to{\rm B}'_{\Pi}$$
	telle que $\omega'\circ\Phi\cong\omega$.
\end{subprop}
\dem Contentons-nous de d\'{e}crire $\Phi$ et un quasi-inverse
$\Psi:{\rm
B}'_{\Pi}\to{\rm B}_{\Pi}$ de $\Phi$.

Soit ${\cal X}=({\sf G',X,G})$ un objet de ${\rm B}_{\Pi}$, et
posons $(G', X,G)=\omega({\cal X})=({\sf |G'|,|X|,|G|})$. Comme le
$\Pi$-groupe
$\sf G$ est constant,
il est donn\'{e} par l'action triviale de $\Pi$ sur le groupe $G$.
Ceci entra\^{\i}ne que l'action de $\Pi$ sur $X$ \emph{commute\/}
\`{a}
celle de $G$, et est donc donn\'{e}e par un morphisme de groupes
$\theta$ de $\Pi$ vers le groupe ${\rm Aut}_{G}(X)$, qui n'est autre
que $G'$, \`{a} isomorphisme canonique pr\`{e}s. On obtient bien
ainsi un objet
$\Phi({\cal X})=(G',X,G,\theta)$ de ${\rm B}'_{\Pi}$.

Inversement, soit $(G',X,G,\theta)$ un objet de ${\rm B}'_{\Pi}$. On
en d\'{e}duit
une action \`{a} gauche de $\Pi$ sur $X$, via $\theta$ et l'action de
$G'$, et une action \`{a} gauche de $\Pi$ sur $G'$ par automorphismes
int\'{e}rieurs via $\theta$. D'o\`{u} un $\Pi$-groupe $\sf G'$ (qui
ne d\'{e}pend
d'ailleurs pas de $X$) et un $\sf G'$-torseur \`{a} gauche $\sf X$.
On
v\'{e}rifie alors imm\'{e}diatement que l'action de $\Pi$ sur $G$
(identifi\'{e} \`{a} ${\rm Aut}_{{G'}}({X})$) est triviale, de sorte
que $\sf
X$ est bien muni d'une structure de $({\sf
G'},\ssoul{G})$-bitorseur.\qed

\begin{subrem}\label{RemConst2}
	Lorsque l'on adopte les descriptions de ${\rm B}_{\Pi}$ et ${\rm
	B}'_{\Pi}$ donn\'{e}es en \ref{RemConst1}, en remarquant en outre
	qu'avec les notations habituelles, on a un isomorphisme $G\cong G'$
	bien d\'{e}fini \`{a} conjugaison pr\`{e}s, on retrouve la bijection
bien
	connue entre ${\rm H}^1(\Pi,G)$ et le quotient de ${\rm
Hom\,}(\Pi,G)$
	par les automorphismes int\'{e}rieurs de $G$.
\end{subrem}
\Subsubsection{Propri\'{e}t\'{e}s de l'\'{e}quivalence de
\ref{EquivBitorsConst}: connexit\'{e}. }\label{EquivConn}
Rappelons qu'un objet $\sf X$ du topos
${\cal C}_{\Pi}$ est \emph{connexe\/} si et seulement si c'est un
$\Pi$-ensemble (non vide et) \emph{transitif}.

Soient ${\cal X}=({\sf G',X,G})$ un objet de ${\rm B}_{\Pi}$, et
$\Phi({\cal X})=(G',X,G,\theta)$ l'objet  de ${\rm B}'_{\Pi}$
correspondant. Nous dirons que ${\cal X}$ est \emph{connexe\/} si
l'objet $\sf X$ de ${\cal C}_{\Pi}$ l'est: cette condition
\'{e}quivaut \`{a} dire que $\theta:\Pi\to G'$ est \emph{surjectif}.

De plus, il existe toujours un objet connexe ${\cal Y}$ de ${\rm
B}_{\Pi}$
(d'ailleurs unique \`{a} isomorphisme non unique pr\`{e}s) et un
morphisme injectif $\cal
Y\to X$. En effet, soit $H'\subset G'$ l'image
de $\theta:\Pi\to G'$; le choix d'un \'{e}l\'{e}ment de $X$
d\'{e}termine un
sous-groupe correspondant $H$ de $G$, et un morphisme
$(H',Y,H)\to(G',X,G)$ de torseurs, qui de plus est un morphisme
$(H',Y,H,\theta_{H'})\to(G',X,G,\theta)$ d'objets de ${\rm B}'_{\Pi}$
o\`{u} $\theta_{H'}$ est obtenu \`{a} partir de $\theta$ par
restriction du but.
Appliquant le foncteur $\Psi$ on obtient le morphisme
$\sf(H',Y,H)\to({G',X,G})$ cherch\'{e}.

\Subsubsection{Propri\'{e}t\'{e}s de l'\'{e}quivalence de
\ref{EquivBitorsConst}: sous-groupes distingu\'{e}s.
}\label{EquivDist}
Soient ${\cal X}=({\sf G',X,G})$  et
$\Phi({\cal X})=(G',X,G,\theta)$ comme en \ref{EquivConn}. Il est
clair que tout sous-groupe distingu\'{e} de $G$ (resp. de $G'$) est
invariant par l'action de $\Pi$, et d\'{e}finit donc un sous-groupe
distingu\'{e} de $\sf G$ (resp. de $\sf G'$). Compte tenu de
\ref{SubSecTriv}, on peut donc identifier canoniquement les quatre
ensembles de sous-groupes distingu\'{e}s de $G$, $G'$, $\sf G$ et
$\sf G'$. De
plus cette identification ne d\'{e}pend pas de $\theta$.

\section{D\'{e}vissage de bitorseurs; preuve du th\'{e}or\`{e}me
\ref{ThPpal}}
\label{SecDem}
\Subsection{Notations et hypoth\`{e}ses. }
\label{SsecForNot}

On reprend les notations et conventions de \ref{NotPiEns} sur les
$\Pi$-ensembles, et l'on se donne une suite exacte
de groupes profinis
\begin{equation}\label{DevPiTer}
	1\ffl \Gamma \ffl \Pi \ffl \pi \ffl 1.
\end{equation}
Nous identifierons la cat\'{e}gorie $\mathcal{C}_{\pi}$ des
$\pi$-ensembles \`{a} la
sous-cat\'{e}gorie pleine de  $\mathcal{C}_{\Pi}$ form\'{e}e des
$\Pi$-ensembles sur lesquels $\Gamma$ op\`{e}re trivialement.

\begin{defi}\label{DefDecomp}
	Si $\mathcal{X}=({\sf G',X,G})$ est un bitorseur dans
$\mathcal{C}_{\Pi}$, nous dirons pour abr\'{e}ger que
\begin{romlist}
	\item\label{DefDecomp1} $\mathcal{X}$ est \emph{de type $\pi$} s'il
	provient d'un bitorseur de $\mathcal{C}_{\pi}$;
	\item\label{DefDecomp2} $\mathcal{X}$ est \emph{de type $\Gamma$}
	s'il existe un morphisme
	$\Phi:\mathcal{X}_{1}=({\sf H',X_{1},H})\to\mathcal{X}$ o\`{u}
	$\sf H'$ est un quotient de $\Gamma$;
	\item\label{DefDecomp3} $\mathcal{X}$ est \emph{d\'{e}composable}
s'il
	peut s'\'{e}crire, \`{a} isomorphisme pr\`{e}s,
		$$\mathcal{X}=\mathcal{Y}\wedge\mathcal{Z}$$
	o\`{u} $\mathcal{Y}$ est de type $\Gamma$ et o\`{u} $\mathcal{Z}$
est de
	type
	$\pi$.
\end{romlist}
\end{defi}

\begin{rem}\label{RemDecomp} Gardons les notations de \ref{DefDecomp}.
	\begin{romlist}
		\item\label{RemDecomp1} Pour que $\mathcal{X}$ soit 
de type $\pi$,
il
		faut et il suffit que $\Gamma$ op\`{e}re trivialement 
sur $\vert\sf
X\vert$.
		\item\label{RemDecomp2} Si $\mathcal{X}$ est de type 
$\Gamma$, le
		morphisme $\Phi$ de \rref{DefDecomp}\,\rref{DefDecomp2} peut
\^{e}tre choisi
		\emph{injectif} (consid\'{e}rer son image, d\'{e}finie en
\ref{MorBitors}).
		\item\label{RemDecomp3} Si $\mathcal{X}$ est d\'{e}composable,
alors
		$\sf G$ est n\'{e}cessairement un $\pi$-groupe.
	\end{romlist}
\end{rem}

\begin{lem}\label{LemDecomp}Soit
$\Phi=(\varphi',u,\varphi):\mathcal{X}=({\sf G',X,G})\to
	\mathcal{X}_{1}=({\sf G_{1}',X_{1},G_{1}})$
	un morphisme de  bitorseurs de $\mathcal{C}_{\Pi}$.
	\begin{romlist}
		\item\label{LemDecomp1} Si $\mathcal{X}$ est de type 
$\Gamma$, il
		en est de m\^{e}me de $\mathcal{X}_{1}$.
		\item\label{LemDecomp2} Si $\mathcal{X}$ est de type 
$\pi$, alors
		$\mathcal{X}_{1}$ est de type $\pi$ si et seulement 
si $\sf G_{1}$
		est un $\pi$-groupe.
		\item\label{LemDecomp3} Si $\mathcal{X}$ est d\'{e}composable,
alors
		$\mathcal{X}_{1}$ est d\'{e}composable si et seulement si $\sf
G_{1}$
		est un $\pi$-groupe.
	\end{romlist}
\end{lem}
\dem L'assertion \ref{LemDecomp1} r\'{e}sulte trivialement de la
d\'{e}finition.
Dans \ref{LemDecomp2} et \ref{LemDecomp3}, le \og seulement si\fg\
est
trivial. R\'{e}ciproquement, pour \ref{LemDecomp2},
remarquer que $\mathcal{X}_{1}$ s'identifie \`{a}
$\mathcal{X}^\varphi$.
La partie \og si\fg\ de \ref{LemDecomp3} r\'{e}sulte de
\ref{LemDecomp1}
et \ref{LemDecomp2} et du lemme \ref{MorProd}.\qed

\begin{thm}\label{ThDeviss}
	Avec les hypoth\`{e}ses et notations de \rref{SsecForNot}, on suppose
	de plus que la suite exacte {\rm(\ref{DevPiTer})} est 
\emph{scind\'{e}e}.

	Alors,
	pour tout groupe fini $G$, tout $(\ast,\ssoul{G})$-bitorseur de
	$\mathcal{C}_{\Pi}$ est d\'{e}composable.
\end{thm}
\dem Soient $G$ un groupe fini et ${\cal X}=({\sf G',X},\ssoul{G})$ un
$(\ast,\ssoul{G})$-bitorseur.
C'est un objet de la cat\'{e}gorie ${\rm B}_\Pi$ de \ref{SorConst}.
Par \ref{EquivConn}, il existe un morphisme
$\mathcal{Y}\to\mathcal{X}$ o\`{u} $\mathcal{Y}$ est un objet connexe de
${\rm B}_\Pi$; par \ref{LemDecomp}\,\ref{LemDecomp3},
$\mathcal{X}$ est d\'{e}composable si $\mathcal{Y}$ l'est.
On peut donc supposer $\mathcal{X}$ \emph{connexe}. Il lui
correspond par \ref{EquivBitorsConst} un objet $(G',X,G,\theta)$
de ${\rm B}'_\Pi$: rappelons que $(G',X,G)$ est le bitorseur ensembliste
$|\cal X|$ et que $\theta:\Pi\to G'$ est un morphisme continu,
ici
\emph{surjectif} puisque $\mathcal{X}$ est connexe.

On d\'{e}duit alors de (\ref{DevPiTer}) et de $\theta$ un diagramme
de $\Pi$-groupes profinis, \`{a} carr\'{e}s commutatifs et \`{a}
lignes exactes:
\setbox0=\hbox{$\ffl$}
\begin{equation}\label{DevG'}
	\begin{array}{rcccccccl}
		1 & \ffl & \Gamma  & \ffl  & \Pi  & \sffl{p}  & \pi & 
\ffl &  1\cr
		& &\mapdown{} & & \mapdown{\theta} &
		\unitlength=\wd0
		\begin{picture}(1,.7)(0,.3)
			\put(.9,.9){\vector(-1,-1){.8}}
			\put(0.3,0.7){\makebox(0,0){$\scriptstyle s$}}
		\end{picture}
		& \mapdown{\overline{\theta}} & & \cr
		1 & \ffl & H' & \ffl & G' & \sffl{q} & \overline{G'} & \ffl & 1
	\end{array}
\end{equation}
o\`{u} les fl\`{e}ches verticales sont surjectives, et o\`{u} $s$
(d\'{e}duit
d'un scindage de (\ref{DevPiTer})) v\'{e}rifie $q\circ
s=\overline{\theta}$.
De plus, le sous-groupe distingu\'{e} $H'$ de $G'$ donne naissance
(par \ref{ClasConj})
\`{a} un sous-$\Pi$-groupe $\sf H'$ de $\sf G'$ et \`{a} un
sous-groupe
distingu\'{e} $H$ de $G$. Posons alors
$$\widetilde{\theta}=s\circ p:\Pi\ffl G';$$
on en d\'{e}duit, par le foncteur $\Psi$ de \ref{EquivBitorsConst},
un bitorseur
$$
{\cal Z}=\Psi(G',X,G,\widetilde{\theta})=({\sf
G''},{\sf Z},\ssoul{G})
$$
avec $|{\cal Z}|=|{\cal X}|$.
Par construction, $\Gamma\subset\ker\widetilde{\theta}$, de sorte
que
${\cal Z}$ est de type $\pi$, et il suffit pour conclure
de voir que le $(\sf G',G'')$-bitorseur
$$
{\cal Y}:={\cal X}\wedge{\cal Z}^{-1}=\soul{\rm Isom}_{G}({\cal
Z},{\cal X})
$$
est de type $\Gamma$. Comme on a
$q\circ\widetilde{\theta}=q\circ{\theta}$,
il est clair que les bitorseurs ${\sf H'}\backslash{\cal X}={\cal
X}/\ssoul{H}$ et ${\cal Z}/\ssoul{H}$ sont isomorphes. Par
suite, d'apr\`{e}s \ref{LemIsomQuot}, il existe un $({\sf
H'},\ast)$-bitorseur $\cal W$ et un morphisme $\cal W\to Y$.
Comme ${\sf H'}$ est  quotient de $\Gamma$, le th\'{e}or\`{e}me est
d\'{e}montr\'{e}.\qed

\Subsection{D\'{e}monstration du th\'{e}or\`{e}me
\ref{ThPpal}.}\label{DemThPpal}
Sous les hypoth\`{e}ses du th\'{e}or\`{e}me \ref{ThPpal}, on
identifie $\mathcal{C}_{\Pi}$  \`{a} la
cat\'{e}gorie des $K$-sch\'{e}mas finis \'{e}tales trivialis\'{e}s
par $M$. Notons $\gre$, comme en
\ref{PropREqBit}\,\ref{PropREqBit2}, la classe des bitorseurs
$\sf(G',X,G)$ de $\mathcal{C}_{\Pi}$ tels que le $\sf G$-torseur $\sf
X$ soit R-\'{e}quivalent au torseur trivial.

Il s'agit de montrer que tout
$(\ast,\ssoul{G})$-bitorseur est dans $\gre$. Soit donc
$\mathcal{X}=({\sf G',X,}\ssoul{G})$ un tel bitorseur. D'apr\`{e}s
\ref{ThDeviss}, on a $\mathcal{X}\cong\mathcal{Y}\wedge\mathcal{Z}$
avec $\mathcal{Y}={\sf(G',Y,G'')}$ de type $\Gamma$ et
$\mathcal{Z}=({\sf G'',Z},\ssoul{G})$ de type $\pi$. Dans ces
conditions, $\sf Z$ est un $\pi$-ensemble, donc $\mathcal{Z}\in\gre$
d'apr\`{e}s la condition \ref{ThPpal1} de \ref{ThPpal}. D'autre part,
$\sf Y$ est induit, comme $\sf G'$-torseur \`{a} gauche, par un
torseur
sous un quotient $\sf H'$ de $\Gamma$, dont le groupe sous-jacent est
un sous-groupe de $|\sf G'|$ (d'apr\`{e}s la remarque
\ref{RemDecomp}\,\ref{RemDecomp2}), donc est isomorphe \`{a} un
sous-groupe de $G$ puisque $|\sf G'|$ est isomorphe \`{a} $G$. Donc
$\mathcal{Y}\in\gre$, vu l'hypoth\`{e}se \ref{ThPpal2} de \ref{ThPpal},
d'o\`{u} aussi $\mathcal{X}\in\gre$ puisque $\gre$ est stable par
produit
(\ref{PropREqBit}\,\ref{PropREqBit2}), et le th\'{e}or\`{e}me est
d\'{e}montr\'{e}.\qed
\bigskip

\noindent{\slshape L'auteur remercie Jean-Louis 
Colliot-Th\'{e}l\`{e}ne et Philippe
Gille pour leurs remarques.}

\end{document}